\documentclass{svproc}
\usepackage{url}

\usepackage{fancyhdr}
\usepackage{bm}
\usepackage{cancel}
\usepackage{latexsym}
\usepackage{amssymb,amsmath,amsfonts,bbm,pifont,upgreek,bbold%,accents
}  %Luigi
\usepackage[colorlinks=true]{hyperref}
\hypersetup{urlcolor=blue, citecolor=red}

\usepackage[dvips]{graphicx}
\newcommand\beq[1]{\begin{equation}\label{#1} }
\newcommand{\eeq}{\end{equation} }

\newcommand\beqa[1]{\begin{eqnarray} \label{#1}}
\newcommand{\eeqa}{\end{eqnarray} }
\newcommand{\beqano}{\begin{eqnarray*} }
\newcommand{\eeqano}{\end{eqnarray*} }
\newcommand\arr[1]{\left\{\begin{array}{l}#1\end{array}\right.}

\newcommand\equ[1]{{\rm (\ref{#1})}}

\expandafter\chardef\csname pre amssym.def
at\endcsname=\the\catcode`\@
\catcode`\@=11
\def\undefine#1{\let#1\undefined}
\def\newsymbol#1#2#3#4#5{\let\next@\relax
 \ifnum#2=\@ne\let\next@\msafam@\else
 \ifnum#2=\tw@\let\next@\msbfam@\fi\fi
 \mathchardef#1="#3\next@#4#5}
\def\mathhexbox@#1#2#3{\relax
 \ifmmode\mathpalette{}{\m@th\mathchar"#1#2#3}%
 \else\leavevmode\hbox{$\m@th\mathchar"#1#2#3$}\fi}
\def\hexnumber@#1{\ifcase#1 0\or 1\or 2\or 3\or 4\or 5\or 6\or 7\or
8\or
 9\or A\or B\or C\or D\or E\or F\fi}
\ifcase\@ptsize
 \font\tenmsb=msbm10
 \font\sevenmsb=msbm7
 \font\fivemsb=msbm5
\or
 \font\tenmsb=msbm10 scaled \magstephalf
 \font\sevenmsb=msbm7 scaled \magstephalf
 \font\fivemsb=msbm5  scaled \magstephalf
\or
 \font\tenmsb=msbm10 scaled \magstep1
 \font\sevenmsb=msbm7 scaled \magstep1
 \font\fivemsb=msbm5 scaled \magstep1
\fi
\newfam\msbfam
\textfont\msbfam=\tenmsb
\scriptfont\msbfam=\sevenmsb
\scriptscriptfont\msbfam=\fivemsb
\edef\msbfam@{\hexnumber@\msbfam}
\def\Bbb#1{\fam\msbfam\relax#1}
\def\widehat#1{\setboxz@h{$\m@th#1$}%
 \ifdim\wdz@>\tw@ em\mathaccent"0\msbfam@5B{#1}%
 \else\mathaccent"0362{#1}\fi}
\def\widetilde#1{\setboxz@h{$\m@th#1$}%
 \ifdim\wdz@>\tw@ em\mathaccent"0\msbfam@5D{#1}%
 \else\mathaccent"0365{#1}\fi}

\def\RIfM@{\relax\ifmmode}
\def\nonmatherr@#1{\errmessage{\string#1\space allowed only in math mode}}
\def\Bbb{\RIfM@\expandafter\Bbb@\else
 \expandafter\nonmatherr@\expandafter\Bbb\fi}
\def\Bbb@#1{{\Bbb@@{#1}}}
\def\Bbb@@#1{\fam\msbfam\relax#1}
\def\setboxz@h{\setbox\z@\hbox}
\def\wdz@{\wd\z@}
\catcode`\@=\csname pre amssym.def at\endcsname
\newcommand{\ii}{{\rm i}  }
\newcommand{\ie}{{\rm i.e.\,}}

\newcommand{\nl}{{\smallskip\noindent}}

\newcommand{\dst}{\displaystyle}
\newcommand\ovl[1]{\overline {#1} }

\newcommand\su[1]{\frac{1}{{#1}} }

\newcommand{\torus}{{\Bbb T}   }
\renewcommand{\natural}{{\Bbb N}   }
\newcommand{\real}{{\Bbb R}   }

\newcommand{\complex}{{\Bbb C}   }

\renewcommand{\a }{{\alpha}   }
\renewcommand{\b}{{\beta}   }
\newcommand{\g}{{\gamma}   }
\newcommand{\G}{{\Gamma}   }
\renewcommand{\d}{{\delta}   }
\newcommand{\D}{{\Delta}   }

\renewcommand{\k}{{\kappa}   }
\renewcommand{\l}{{\lambda}   }
\renewcommand{\L}{{\Lambda}   }
\newcommand{\m}{{\mu}   }
\newcommand{\n}{{\nu}   }
\newcommand{\x }{{\xi}   }

\newcommand{\p}{{\pi}   }
\renewcommand{\P}{{\Pi}   }
\renewcommand{\r}{{\rho}   }
\newcommand{\s}{{\sigma}   }

\renewcommand{\t}{{\tau}   }
\newcommand{\f}{{\varphi}   }

\renewcommand{\o}{{\omega}   }

\renewcommand{\Im}{{\, \rm Im\, }}

\newcommand{\cB}{{\cal B} }
\newcommand{{\cE}}{{\cal  E} }
\newcommand{\cT}{{\cal T} }

\newcommand{{\cH}}{{\cal H} }
\newcommand{{\cK}}{{\cal K} }
\newcommand{\cC}{{\cal C} }
\newcommand{\cD}{{\cal D} }
\newcommand{\cF }{{\cal F} }
\newcommand{\cG}{{\cal G} }
\newcommand{{\cJ}}{{\cal J}}
\newcommand{\cL}{{\cal L} }
\newcommand{\cM}{{\cal M} }
\newcommand{\cP}{{\cal P} }

\newcommand{\cN}{{\cal N} }
\newcommand{\cS}{{\cal S} }

\newcommand\bx{{\mathbf x}}

\newcommand\by{{\mathbf y}}

\newcommand\ppu{{(1) }}
\newcommand\ppd{{(2) }}

\newcommand{\CC}{{\rm C}}

\newcommand{\EE}{{\rm E}}
\newcommand\FF{{\rm F}}
\newcommand\GG{{\rm G}}
\newcommand\HH{{\rm H}}

\newcommand\JJ{{\rm J}}

\newcommand\RR{{\rm R}}

\newcommand\UU{{\rm U}}
\newcommand\VV{{\rm V}}

\newcommand\ZZ{{\rm Z}}
\newcommand\dd{{\rm d}}
\newcommand\ee{{\rm e}}
\newcommand\hh{{\rm h}}

\newcommand\rr{{\rm r}}
\usepackage{color}
\definecolor{yellow}{rgb}{0.99, 0.93, 0.0}

\newcommand{\R}{\mathbb{R}}

\newcommand\mm{{\rm m}}

\newcommand\cO{{\cal O}}

\begin{document}
\mainmatter              % start of a contribution
\title{A new analysis of the three--body problem}
\titlerunning{Three--body problem}  % abbreviated title (for running head)
\author{J\'er\^ome Daquin\inst{1}
 \and Sara Di Ruzza\inst{2} \and  Gabriella Pinzari\inst{2}}
\authorrunning{J. Daquin, S. Di Ruzza,  and G. Pinzari} % abbreviated author list (for running head)
\tocauthor{J\'er\^ome Daquin, Sara Di Ruzza,  and Gabriella Pinzari}
\institute{ %WWW home page:
\texttt{%http://users/\homedir iekeland/web/welcome.html
}
naXys Research Institute, University of Namur, Belgium
\email{jerome.daquin@unamur.be} \and Department of Mathematics, University of Padua, Italy\\
\email{sara.diruzza@unipd.it, gabriella.pinzari@math.unipd.it}}
\maketitle              % typeset the title of the contribution

\begin{abstract}
In the recent papers~\cite{pinzari20a},~\cite{diruzzaDP20}, respectively, the existence of motions where the perihelions afford periodic oscillations about certain equilibria and the onset of a topological horseshoe have been proved. Such results have been obtained using, as neighbouring integrable system, the so--called two--centre (or {\it Euler}) problem and a suitable canonical setting proposed in~\cite{pinzari19},~\cite{pinzari20b}. Here we review such results.
% We would like to encourage you to list your keywords within
% the abstract section using the \keywords{...} command.
\end{abstract}
\section{Overview}\label{2centers}
{%\bf\large 1.2} {
In the recent papers~\cite{diruzzaDP20},~\cite{pinzari20a} the existence, in the three--body problem (3BP), of motions which by no means  can be regarded as ``extending'' in some way Keplerian motions has been proved. Indeed, the motions found in those papers can be  better understood as continuations of the motions of the so--called {\it two--centre problem} (or {\it Euler  problem}; 2CP from now on).\\
The motivation that pushed such researches was a new analysis of 2CP  carried out in~\cite{pinzari20b}, combined with a remarkable property -- which we called {\it renormalizable integrability} -- pointed out in~\cite{pinzari19}. It relates the ``simply averaged  Newtonian potential'' (see the precise definition below) and the function, which in this paper we shall refer to as {\it Euler integral}, that makes the 2CP integrable. Roughly,  such property states that the averaged Newtonian potential and the Euler integral have the same motions, as they are one a function of the other. As, on the other hand, the motions of the Euler integral are, at least qualitatively, explicit, and the averaged Newtonian potential is  a prominent  part  of the 3BP Hamiltonian,  the papers~\cite{diruzzaDP20},~\cite{pinzari20a} gave partial answers to the natural question whether the motions of the Euler integral can be traced in 3BP. Let us introduce some mathematical tools in order to make our statements more precise.

\nl
In terms of Jacobi coordinates~\cite{giorgilli}    the three--body problem Hamiltonian with masses $1$, $\m$, $\k$ is the translation--free function}
\beqano
{\HH_{\rm J}}&=&\frac{\|\mathbf y\|^2}{2}\left(1+\frac{1}{\m}\right)+\frac{\|\mathbf y'\|^2}{2}\left(\frac{1}{1+\m}+\frac{1}{\k}\right)-\frac{ \m }{\|\mathbf x\|}-\frac{\m \k }{\|\mathbf x'-\frac{1}{1+\m}\mathbf x\|}\nonumber\\
&&-\frac{ \k }{\|\mathbf x'+\frac{\m}{1+\m}\mathbf x\|}\, .
\eeqano
{Here, $(\by', \by, \bx', \bx)\in (\real^3)^4$ (or $(\real^2)^4$, in the planar case), $\|\cdot\|$ denotes Euclidean norm and the gravity constant has been taken equal to one, by a proper choice of the units system.}
{We rescale impulses and positions}
\beqa{resc}{\mathbf y\to \frac{\m}{1+\m}\mathbf y\, , \quad \mathbf x\to {(1+\m)}{\mathbf x}\, , \quad \mathbf y'\to \m\b{\mathbf y'}\, , \quad \mathbf x'\to \b^{-1}\mathbf x'},\eeqa
{multiply the Hamiltonian by $\frac{1+\m}{\m}$ (by a rescaling of time) and obtain
\beqa{Jac}{\HH_{\rm J}=\frac{\|\mathbf y\|^2}{2}-\frac{1}{\|\mathbf x\|}+\g \left(\frac{\|\mathbf y'\|^2}{2}
-\frac{\ovl\b}{\b+\ovl\b}\frac{ 1 }{\|\mathbf x'-\b\mathbf x\|}
-\frac{\b}{\b+\ovl\b}\frac{ 1 }{\|\mathbf x'+\ovl\b\mathbf x\|}
\right)},\eeqa
 with}
\beqa{betagamma1}{\g =\frac{\k^3(1+\m)^4}{\m^3(1+\m+\k)}\, , \quad 
\b=\frac{\k^2(1+\m)^2}{\m^2(1+\m+\k)}\, , \quad \ovl\b=\m\b}%\, , \quad \ovl\m:=\frac{1}{1+\m}
\, . \eeqa
Likewise, one might consider 
the problem written in the so--called $1$--centric coordinates. In that case,
\beqa{helio}\HH_{0}&=&\frac{\|{\mathbf y}\|^2}{2}-\frac{1}{\|{\mathbf x}\|}+\g\left(
\frac{\|{\mathbf y}'\|^2}{2}-\frac{\b}{\b+\ovl\b}\frac{ 1}{\|{\mathbf x}'\|}-\frac{\ovl\b}{\b+\ovl\b}\frac{1}{\|{\mathbf x}'-(\b+\ovl\b){\mathbf x}\|}
\right)\nonumber\\
&&+\ovl\b\, {{\mathbf y}'\cdot {\mathbf y}},\eeqa
with $\g$, $\b$ and $\ovl\b$ analogous to~\equ{betagamma1}, up to replace the  factors $(1+\m+\k)$  with $(1+\k)$.
Note that we are not assuming $\m$, $\k\ll1$ (in fact, in our applications, we shall make different choices), which means that Jacobi or $1$--centric coordinates above are not necessarily centered at the most massive body.  In order to simplify the analysis, we introduce a main assumption.
Both the Hamiltonians $\HH_{\rm J}$ and $\HH_0$ in~\equ{Jac} and~\equ{helio} include the Keplerian term  
\beqa{Kep}\JJ_0:=\frac{\|\mathbf y\|^2}{2}-\frac{1}{\|\mathbf x\|}=-\frac{1}{2\L^2}\, . \eeqa
We  assume that  $\JJ_0$ is a ``leading'' term in such  Hamiltonians. By averaging theory, this assumption allows us to replace
(at the cost of a small error) $\HH_{\rm J}$ and $\HH_0$ with  their respective  $\ell$--averages
\beqa{ovlH}\ovl\HH_{i}=-\frac{1}{2\L^2}+\g\widehat\HH_{i}\, \eeqa
 {with $i$=J,0,} where $\ell$ is the mean anomaly associated to~\equ{Kep}, and\footnote{Remark that 
  ${\mathbf y}(\ell)$ has  vanishing $\ell$--average so that the last term in~\equ{helio} does not survive.}
 \beqa{secular}
\widehat\HH_\JJ&:=&\frac{\|{\mathbf y}'\|^2}{2 }-\frac{\ovl\b}{\b+\ovl\b}\UU%^+
_{\b}-\frac{\b}{\b+\ovl\b}\UU%^-
_{-\ovl\b}\nonumber\\
\widehat\HH_0&:=&\frac{\|{\mathbf y}'\|^2}{2 }-\frac{\ovl\b}{\b+\ovl\b}\UU%^+
_{\b+\ovl\b}-\frac{\b}{\b+\ovl\b}\frac{1}{\|{\mathbf x}'\|}
\eeqa
 with
\beqa{Usb}\UU
_{\b}:=\frac{1}{2\p}\int_0^{2\p}\frac{d\ell}{\|{\mathbf x}'-\b{\mathbf x}(\ell)\|}.\eeqa
 {In these formulae, the term $-\su{2\L^2}$ will be referred to as ``Keplerian term'', while terms of the form $-\frac{1}{\|\bx'-\b\bx\|}$ will be called ``Newtonian potentials''. Therefore, $\UU_\b$ will be called ``averaged Newtonian potential''.}    What we want to underline in that respect is that the averages~\equ{ovlH} are ``simple'', i.e., computed with respect to only one  mean anomaly. Most often, in the literature double averages are considered; e.g.~\cite{laskarR95},~\cite{fejoz04},~\cite{pinzari-th09},~\cite{chierchiaPi11b}, ~\cite{fejozG2016},~\cite{gioLS17}.

\nl{Whether and at which extent the Hamiltonians \eqref{ovlH} are good approximations of \eqref{Jac}, \eqref{helio} is a demanding question, as, besides the mass parameters $\m$, $\k$, also the region of phase space which is being considered plays a crucial r\^ole. We limit ourselves to some heuristics, focusing, in particular, on the case considered in~\cite{diruzzaDP20}. Here the Hamiltonian \eqref{Jac} has been investigated, with $\m=1\gg\k$, or, equivalently, $\b=\ovl\b\gg 1$ (see \eqref{Cbetaovlbeta} for the precise values). Physically, this corresponds to a couple of asteroids with equal mass interacting with a star, with $\bx$ being the relative distance of the twin asteroids, and $\bx'$ the distance of the star from their center of mass. In that case, the region of phase space was chosen so that $\|\bx'\|>\b\|\bx\|$, so that the two denominators of the Newtonian potentials do not vanish. Expanding such Newtonian potentials  in powers of $\frac{\b a}{\rr}$, where $a=\L^2$, $\rr:=\|\bx'\|$, one sees that the lowest order terms depending on $\ell$ have size $\frac{\gamma\b a}{\rr^2}\sim \frac{\k^3 a}{\rr^2}$ (as $\b\sim \k$, $\g\sim \k^2$). So, such terms are negligible compared to the size $\frac{1}{a}$ of the Keplerian term, provided that $\frac{\k^{3/2}a}{\rr}\ll1$.}

\nl
{We now turn to describe the main features of the Hamiltonians \eqref{ovlH}. } Neglecting the {Keplerian term}, which is an inessential additive constant for $\ovl\HH_i$ and reabsorbing the constant $\g$ with a  time change, we are led  to look at the Hamiltonians $\widehat\HH_i$ in~\equ{secular}, which, from now on, will be our object of study. Without loss\footnote{We can do this as  the Hamiltonians $\HH_\JJ$  and $\HH_0$ rescale by a factor $\b^{-2}$ as $(\by', \by)\to\b^{-1}(\by', \by)$ and $(\bx', \bx)\to\b^2(\bx', \bx)$.} of generality,   we    fix the constant action $\L$ to $1$. 

\nl   
 For definiteness and simplicity, we describe the setting in the case of the planar problem, in which case, after  
 reducing the SO(2) symmetry, $\widehat\HH_i$ have 2 degrees--of--freedom; all the generalisations to the spatial problem  being described in Section~\ref{The classical integration of the two-centre problem}. To describe the coordinates we used, we denote as $\mathbb E$ the Keplerian ellipse generated by  Hamiltonian~\equ{Kep}, for  negative values of the energy. Assume $\mathbb E$ is not a circle. Remark that, as  the mean anomaly  $\ell$ is averaged out, we loose any information concerning the position of $\bx$ on $\mathbb E$, so we shall only need two couples of coordinates for determining the shape of 
 $\mathbb E$ and the vectors $\by'$, $\bx'$. These are:
 
\begin{itemize}
\item[\tiny\textbullet] the  ``Delaunay couple'' $(\GG, {\rm g})$, where $\GG$ is the Euclidean length of $\bx\times \by$ and ${\rm g}$ detects the perihelion. We remark that ${\rm g}$ is measured with respect to $\bx'$ (instead of with respect to a fixed direction), as  the  SO(2) reduction we use 
 fixes a rotating frame which moves with $\bx'$ (compare the formulae in~\equ{coord});

\item[\tiny\textbullet] the ``radial--polar couple''$(\RR, \rr)$, where $\rr:=\|\bx'\|$ and $\RR:=\frac{\by'\cdot\bx'}{\|\bx'\|}$.
%\item[\tiny\textbullet]  $2$  coordinates $(\CC, \g)$ for the total angular momentum of the system. 
\end{itemize}
 %Moreover, due to the conservation of $\CC$ along the trajectories of $\widehat\HH_\JJ$ and $\widehat\HH_0$, $\g$ is cyclic and hence $\CC$ can be regarded as a fixed number.

\begin{figure}
  	\centering
  	\includegraphics[width=0.95\linewidth]{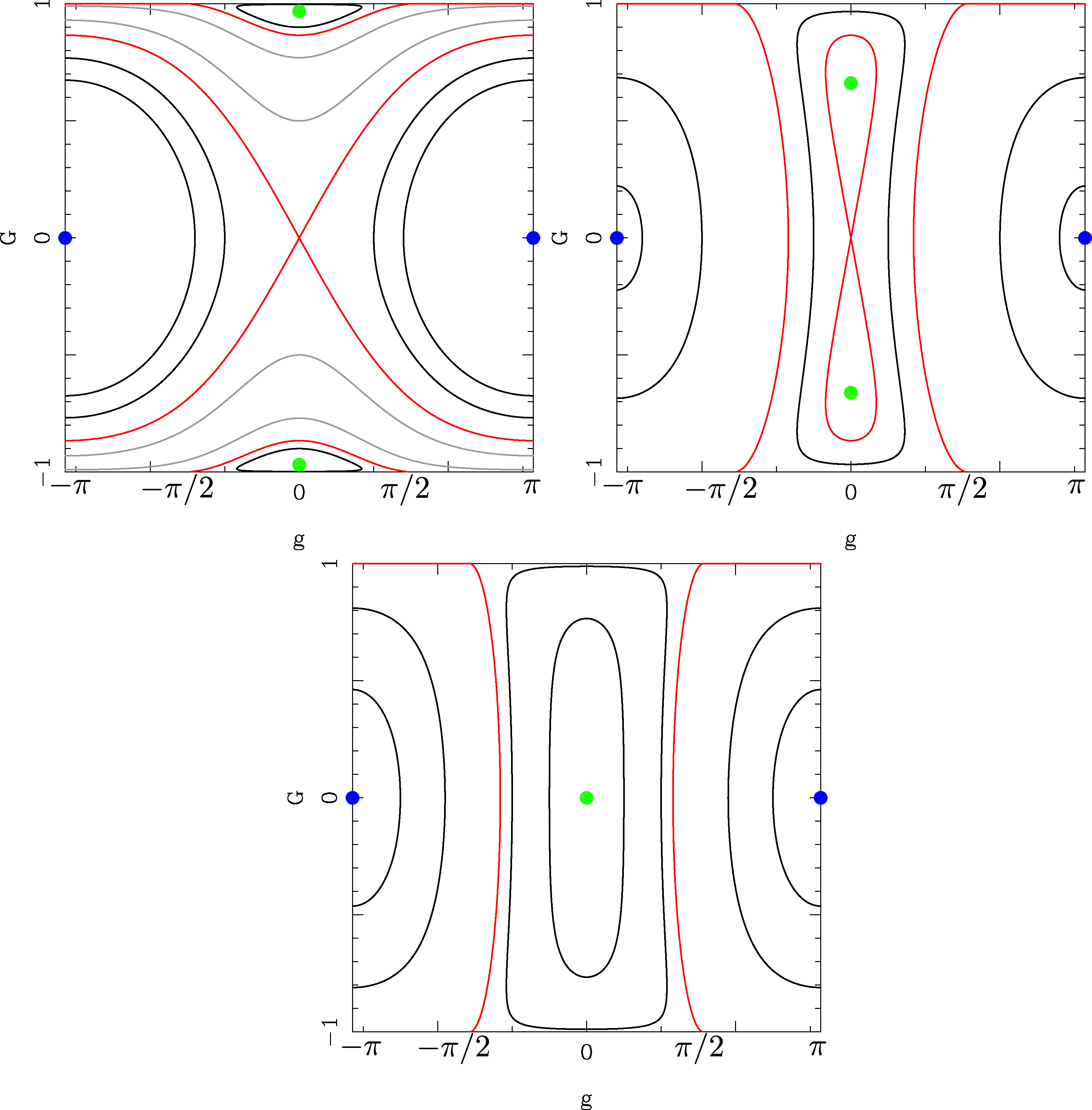}
  	\caption{Phase portraits of $\textrm{E}_{0}$ given by (\ref{eq:E0}) in the plane $(\textrm{g},\textrm{G})$ for 
  		 $0<\textrm{r}<1$ {(top left)}, $1<\textrm{r}<2$ {(top right)} and
  		 $\textrm{r}>2$ {(bottom)}.
  		The points corresponding to minima of $\textrm{E}_{0}$ are labeled in blue, maxima appear in green and separatrices correspond to red curves.} 
  	\label{fig:fig1}
  \end{figure}

\nl
We now describe what we mean by {\it renormalizable integrability}~\cite{pinzari19}.
Note first that, in terms of the coordinates above, the functions $\UU_\b(\rr, \GG, {\rm g})$ in~\equ{Usb} depend on $(\rr, \GG, {\rm g})$ and remark the homogeneity property
\beqa{homogeneity}
\UU_\b(\rr, \GG, {\rm g})=\b^{-1}\UU(\b^{-1}\rr, \GG, {\rm g})\qquad {\rm where}\quad \UU:=\UU_1\, . 
\eeqa
By {\it renormalizable integrability} we mean that  there exists a function $\FF$ of two arguments such that the function $\UU$ in~\equ{homogeneity} verifies
\beqa{relation}\UU(\rr, \GG, {\rm g})=\FF(\EE_0(\rr, \GG, {\rm g}), \rr)\,,\eeqa
where
\beqa{E}\label{eq:E0}
\EE_0(\rr, \GG, {\rm g})=\GG^2+\rr\sqrt{1-\GG^2}\cos{\rm g}\,.\eeqa
 By~\equ{relation},  the level curves of $\EE_0$ are also level curves of $\UU$. On the other hand, the phase portrait of  $\EE_0$ in the plane $({\rm g}, \GG)$ -- i.e., 
the family  of curves
\beqa{level sets}\dst\EE_0(\rr, \GG, {\rm g})=\GG^2+\rr \sqrt{1-\GG^2}\cos{\rm g}=\cE,
\eeqa
in the plane $({\rm g}, \GG)$ accordingly to the different values of $\rr$ -- is completely explicit~\cite{pinzari20b}. For $0<\rr<1$ or $1<\rr<2$ it includes two minima $(\pm\p, 0)$  on the ${\rm g}$--axis; two symmetric maxima on the $\GG$--axis and one saddle point at $(0, 0)$. When $\rr>2$ the saddle point disappears and $(0, 0)$ turns to be a maximum. 
The phase portrait includes two separatrices in the case $0<\rr<2$; one separatrix in the case $\rr>2$. These are the level set $\cS_0(\rr)$ through the saddle, corresponding to $\cE=\rr$, for $0<\rr<2$, and the level set $\cS_1(\rr)=\{\cE=1\}$, for any $\rr$. Rotational motions in between $\cS_0(\rr)$ and $\cS_1(\rr)$, do exist only for $0<\rr<1$. The minima and the maxima are surrounded by librational motions and different motions (librations about different equilibria or rotations) are separated by $\cS_0(\rr)$ and $\cS_1(\rr)$. {The reader is referred to Figure\,\ref{fig:fig1} for further qualitative details about the portion of the phase space corresponding to  $[-\pi,\pi] \times [-1,1]$.}

\nl
 We call  {\it perihelium librations} the librational motions about $(\pm\p, 0)$  or  $(0, 0)$.
  Their physical meaning  is that the  perihelion of ${\mathbb E}$  affords oscillations  while    $\mathbb E$, highly eccentric  anytime, periodically flattens to a segment in correspondence of the times when $\GG$ vanishes. After the flattening time,  the sense of rotation on $\mathbb E$ is reversed (as $\GG$ changes its sign). We remark that (see the next section for a discussion)  the potential $\UU$ is well defined along the level sets of $\EE$, with the exception of $\cS_0(\rr)$, where $\UU$ is singular. In particular, $\UU$ remains regular for all $\rr>2$. 

 \nl
Let $0<\b_*\le \b^*$  be defined via
\beqa{b*}\b_*:=\left\{\begin{array}{lll}
\dst\frac{\b\ovl\b}{\b+\ovl\b}\quad &{\rm for}\ \HH_{\rm J}\\\\
\dst\ovl\b&{\rm for}\ \HH_0
\end{array}
\right.\qquad \b^*:=\left\{\begin{array}{lll}\dst \max\{\b,  \ovl\b\}\quad &{\rm for}\ \HH_{\rm J}\\		\\\\
\dst\b+\ovl\b&{\rm for}\ \HH_0\,.
\end{array}
\right.\eeqa
De--homogeneizating via~\equ{homogeneity}, we see that, if $\rr>2\b^*$, then 
we fall  in the third panel  in Figure~\ref{fig:fig1} for any
 $\UU_\b$'s in~\equ{secular} 
 so that all of such potentials
 afford perihelion librations about $(0, 0)$ and $(\pm \p, 0)$.
The works~\cite{pinzari20a},~\cite{diruzzaDP20} deal precisely with this situation.

\nl
Before (and in order to) describing the purposes of such works, we informally discuss  the r\^ole of the total
 angular momentum's length $\CC:=\|\bx\times\by+\bx'\times\by'\|$.  This quantity enters in~\equ{secular} via kinetic term $\|\by'\|^2$, according to
 \beqa{byp}\|\by'\|^2=\RR^2+\frac{(\CC-\GG)^2}{\rr^2},\eeqa
(as $|\CC-\GG|$ is the Euclidean length of $\bx'\times \by'$, assuming that $\bx\times\by$ and $\bx'\times\by'$ are parallel). Combining~\equ{byp} with 
an expansion 
 \beqano\UU_\b(\rr, \GG, {\rm g})=-\su\rr+\su\rr\sum_{k\ge 1}u_{\k}(\GG, {\rm g})\left(\frac{\b}{\rr}\right)^k,\eeqano
 of the $\UU_\b$'s
 in~\equ{secular} in powers of $\rr^{-1}$, one can
 split the Hamiltonians $\widehat\HH_\JJ$ and $\widehat\HH_0$  in~\equ{secular} in two parts, which we call, respectively, {\it fast} and {\it slow}:
\beqa{split3}\HH_{\rm fast}=\frac{\RR^2}{2}+\frac{\CC^2}{2\rr^2}-\frac{1}{\rr}\,,\quad \HH_{\rm slow}:=\widetilde\UU_{\b, \ovl\b}(\rr, \GG, {\rm g})+\frac{-2\CC\GG+\GG^2}{2\rr^2},\eeqa
where $\widetilde\UU_{\b, \ovl\b}$ collects  terms
of order $\frac{\b^*}{\rr^2}$, or higher, hence, retains the symmetries and the equilibria of the $\UU_\b$'s discussed above. We fix, in phase space, a region of initial data where the terms  in~\equ{split3} verify (see~\cite{diruzzaDP20} for an informal discussion)
\beqa{order}
\|\HH_{\rm fast}\|\gg\|\widetilde\UU_{\b, \ovl\b}\|\gg\|\HH_{\rm slow}-\widetilde\UU_{\b, \ovl\b}\|\,.
\eeqa
Then, the motions of $\RR$ and $\rr$, mainly ruled $\HH_{\rm fast}$, are faster than the ones of $\GG$ and ${\rm g}$, ruled by $\HH_{\rm slow}$. 
 If  $\CC=0$, the smallest term $\HH_{\rm slow}-\widetilde\UU_{\b, \ovl\b}=\frac{\GG^2}{2\rr^2}$ is even with respect to $\GG$, so  $\HH_{\rm slow}$ retains the symmetries and  equilibria of $\widetilde\UU_{\b, \ovl\b}$ in Figure~\ref{fig:fig1}, {for the case $\textrm{r}>2$}. 
However, in this case, 
$\HH_{\rm fast}$ is unbounded below, so
nothing prevents
$\rr$  to decrease below $2\b^*$ and  the  scenario  rapidly change{s}  from  (c) to (b) or (a). In this case, one has then to prove that perihelion librations occur in  the full Hamiltonians~\equ{secular}  in such a short time that it prevents the scenario to change. The following result was obtained:

\begin{theorem}[\cite{pinzari20a}]\label{main}
Take, in~\equ{secular}, $\CC=0$.
Fix an arbitrary neighbourhood $\UU_0$ of $(0,0)$ or of $(0, \p)$ and an arbitrary neighbourhood $\VV_0$ of an unperturbed curve $\g_0(t)=(\GG_0(t), {\rm g}_0(t))\in \UU_0$ in Figure~\ref{fig:fig1}. Then it is possible to find six numbers $0<c<1$, $0<\b_-<\b_+$, $0<\a_-<\a_+$, $T>0$, such that, for any $\b_-<\b_*\le \b^*<\b_+$ 
the projections $\G_0(t)=(\GG(t), {\rm g}(t))$ of
all the orbits $\G(t)=(\RR(t), \GG(t), \rr(t), {\rm g}(t))$ of $\ovl\HH_1$, $\ovl\HH_2$ with initial datum $(\RR_0, \rr_0, \GG_0, {\rm g}_0)\in [\su{\sqrt{c\a_+}}, \su{\sqrt{c\a_-}}]\times [c\a_-, \a_+]\times \UU_0$
belong to  $\VV_0$ for all $0\le t\le T$. Moreover, the angle $\g(t)$  between the position ray of $\G_0(t)$ and the ${\rm g}$--axis affords a variation larger than $2\p$ during the time $T$.
\end{theorem}

\nl
The proof of Theorem~\ref{main} uses a new normal form theorem, together with the construction of a system of coordinates well adapted to perihelion librations, as reviewed in Section~\ref{Asymptotic  action--angle coordinates}.

\nl
If $\CC\ne 0$, $\HH_{\rm fast}$ is bounded below, attaining its minimum at
\beqa{equilibrium}
\RR_0=0\, , \quad \rr_0={\CC^2}\,.
\eeqa
It is reasonable to expect that if the
initial values of $\RR$ and $\rr$ are close to~\equ{equilibrium}, they will remain there for some time and  the motions of
 $\GG$ and ${\rm g}$ will be  close to be ruled by
$\HH_{\rm slow}^0:=\HH_{\rm slow}|_{\rr=\rr_0}$,
{
which reads, to the lowest order{s},
\begin{eqnarray}\label{eq:H0slow}
\textrm{H}_{\textrm{slow}}^{0}
(\textrm{G},\textrm{g})&=&
\frac{{-2\textrm{C}\textrm{G}+\textrm{G}^2}}{2\textrm{r}_{0}^{2}}
{-}
\beta^{2} \frac{
{(5-3\textrm{G}^{2})}
}{{8}\textrm{r}_{0}^{{3}}}
{-}
\beta^{2} \frac{
15{(1-\textrm{G}^{2})}
}{{8}\textrm{r}_{0}^{{3}}}\cos2\textrm{g}\nonumber\\
&&+{\rm O}({\rm r}_0^{-5})\,
.
\end{eqnarray}
}

\begin{figure}
	\centering
	\includegraphics[width=0.5\linewidth]{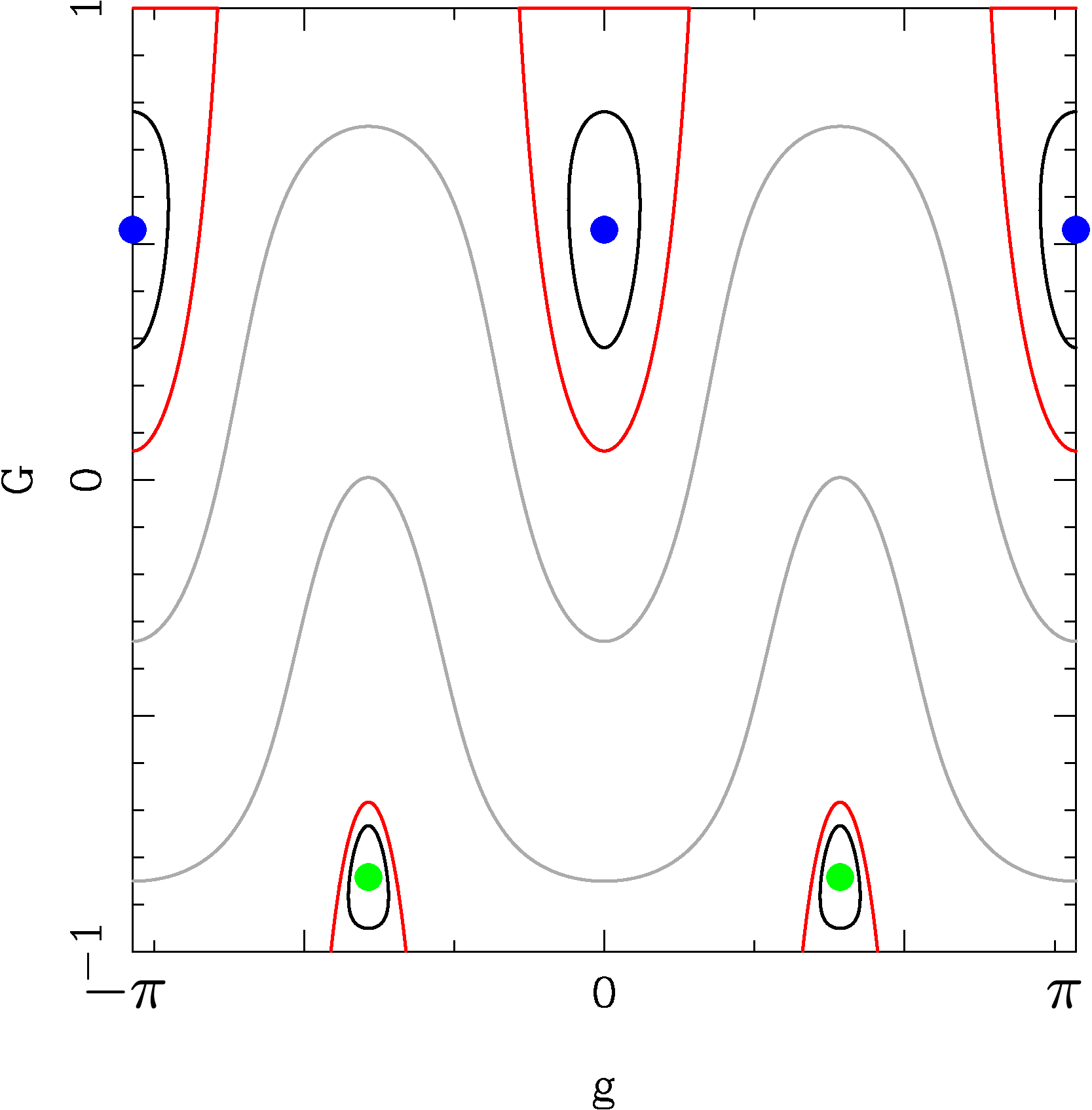}
	\caption{Phase portrait of $\HH_{\rm slow}^0$ given in (\ref{eq:H0slow}) with $C=25, \beta=\bar{\beta}=80$.} 
	\label{fig:fig2}
\end{figure}

On the other hand,
the equilibria of $\widetilde\UU^0_{\b, \ovl\b}$ have some chance of surviving in $\HH^0_{\rm slow}$ for small values of $|\CC|$, but the symmetries of $\widetilde\UU^0_{\b, \ovl\b}$ do not persist.  As an example, in Figure\,~\ref{fig:fig2}, we report the phase portrait  of $\HH_{\rm slow}^0$ for $\widehat\HH_\JJ$, with
\beqa{Cbetaovlbeta}
\CC=25\,,\quad\ \b=\ovl\b=80\,.
\eeqa
We call {\it unperturbed motions}  the motions obtained combining~\equ{equilibrium} with the motions in Figure~\ref{fig:fig2}. The natural question now is 
 whether and at which extent the motions of~\equ{secular} may be regarded as perturbations of such unperturbed ones. 
The question was considered in~\cite{diruzzaDP20}, from the numerical point of view. 
Namely, in~\cite{diruzzaDP20}
the full motions of $\widehat\HH_{\JJ}$ were analysed, with $\CC$, $\b$ and $\ovl\b$   chosen\footnote{\label{comparison}The quantities $\b$, $\ovl\b$, $\g$, $\by'$, $\bx'$, $\by$, $\bx$, $\RR$, $\rr$, $\GG$, ${\rm g}$, $\CC$ of the present paper are related to  $\b$, $\ovl\b$, $\s$, $y'$, $x'$, $y$, $x$, $\RR$, $\rr$, $\GG$, ${\rm g}$, $\CC$
 in~\cite{diruzzaDP20} via the relations (with ``here'' , ``there'' standing for ``in the present paper'' and ``in~\cite{diruzzaDP20}'', respectively) $\b_{\rm here}=(1+\m)\b_{\rm there}$, $\ovl\b_{\rm here}=(1+\m)\ovl\b_{\rm there}$, $\s(1+\m)^2=\g$, $\by'=\frac{\L}{1+\m}y'$, $\bx'=\frac{1+\m}{\L^2}x'$, 
 $\by={\L}y$, $\bx=\frac{x}{\L^2}$,
   $\RR_{\rm here}=\frac{\RR_{\rm there}\L}{(1+\m)}$, $\rr_{\rm here}=\frac{\rr_{\rm there}(1+\m)}{\L^2}$, $\GG_{\rm here}=\frac{\GG_{\rm there}}{\L}$, ${\rm g}_{\rm here}={\rm g}_{\rm there}$, $\CC_{\rm here}=\frac{\CC_{\rm there}}{\L}$ where $\L$, $\m$ were chosen, in~\cite{diruzzaDP20}, $3.099$ and $1$, respectively. 
 Note also the  misprint in the definition of $\s$ in~\cite[(1.4)]{diruzzaDP20}, as the power of $\m$ at the denominator should be $3$ instead of $2$. This misprint is inessential, as the number $\s$ plays no r\^ole in~\cite{diruzzaDP20}.}
 as in~\equ{Cbetaovlbeta}, and the initial values of $\RR$ and $\rr$ close to~\equ{equilibrium}. Numerical evidence of orbits continuing the unperturbed  orbits above, interposed with zones of chaos, was obtained.

  \noindent
 {The paper is organised as follows:
 	\begin{enumerate}
 		\item In section \ref{The classical integration of the two-centre problem}, we review recent results on the two--centre problem. We discuss first the existence of an invariant, referred as \textit{Euler integral}, whose  
 		 expression in the asymmetric setting is also given. Taking advantages of canonical coordinates lowering the number of degrees--of--freedom, coupled together with the \textit{renormalizable integrability} property, the level sets of the averaged Newtonian potential are discussed in the planar case.       
       \item In section \ref{Asymptotic action--angle coordinates}, we outline the proof of Theorem~\ref{main} following~\cite{pinzari20a}. The proof relies on normal form of Hamiltonians~\equ{secular} free of small divisors, combined with an expression of the Euler integral suited for large values of $\rr$.
 		\item In section \ref{Chaos in a binary asteroid system}, we further complement the understanding of the dynamics in the regime of large $\rr$: we retrace the steps of~\cite{diruzzaDP20} in constructing explicitly an horseshoe orbit, therefore introducing the existence of symbolic dynamics. The methodology relies essentially on the construction of ``boxes'' stretching across one another under the action of a specific Poincar\'e mapping, and uses arguments of covering relations as introduced in~\cite{ZGLICZYNSKI200432}.
 	\end{enumerate}
}

 \section{Euler problem revisited}\label{The classical integration of the two-centre problem} 
 
 {In this section we review the classical integration of the two--centre problem and complement it with considerations that will be useful to us in the next.}
  
  		\noindent
 {The} 2CP is the  system, in ${\mathbb R}^3$ (or ${\mathbb R}^2$), of one particle interacting with two fixed masses via Newton Law. If $\pm{\mathbf v}_0\in {\real}^3$ are the position coordinates of the centers, $\mm_\pm$ their masses;
${\mathbf v}$, with  ${\mathbf v}\ne \pm{\mathbf v}_0$, the position coordinate of the moving particle; ${\mathbf u}=\dot{\mathbf v}$ its velocity, and $1$ its mass,
the Hamiltonian of the system ({\it Euler Hamiltonian}) is 
\beqa{2Cold}\JJ=\frac{\|{\mathbf u}\|^2}{2}-\frac{\mm_+}{\|{\mathbf v}+{\mathbf v}_0\|}-\frac{\mm_-}{\|{\mathbf v}-{\mathbf v}_0\|}\, . \eeqa
Euler showed~\cite{jacobi09} that $\JJ$ exhibits  $2$ independent first integrals, in involution. One of these first integrals is the projection 
\beqa{Theta}\Theta= {\mathbf M}\cdot \frac{{\mathbf v}_0}{\|{\mathbf v}_0\|}\eeqa
of the angular momentum ${\mathbf M}={\mathbf v}\times {\mathbf u}$ of the particle along the direction ${\mathbf v}_0$. It 
 is not specifically due to the Newtonian potential, but, rather, to  its invariance  by rotations around the axis ${\mathbf  v}_0$. For example, it  persists if the Newtonian potential is replaced with a $\a$--homogeneous one. 
The existence of the following constant of motion, which we shall refer to as {\it Euler integral}:

\beq{G1}\EE=\|{\mathbf v}\times {\mathbf u}\|^2+({\mathbf v}_0\cdot {\mathbf u})^2+2 {\mathbf v}\cdot {\mathbf v}_0\left(\frac{\mm_+}{\|{\mathbf v}+{\mathbf v}_0\|}-\frac{\mm_-}{\|{\mathbf v}-{\mathbf v}_0\|}\right)\eeq
is pretty specific of $\JJ$. As observed in~\cite{bekovO78}, in the limit of merging centers, i.e., ${\mathbf v}_0=\mathbf 0$, $\JJ$ reduces to the Kepler Hamiltonian~\equ{Kep}, and  $\EE$ to the squared length of the angular momentum of the moving particle.

\nl
The formula in~\equ{G1} is  not easy\footnote{See however~\cite{dullinM16} for a formula related to~\equ{G1}.} to be found in the literature, so we briefly discuss it.\\
After fixing a reference frame with the third axis in the direction of ${\mathbf v}_0$ and denoting as $(v_1,  v_2,  v_3)$ the coordinates of $ {\mathbf v}$ with respect to such frame,
one introduces  
the so--called  ``elliptic coordinates''  

 \beqa{lambdabeta} \l=\frac{1}{2}\left(\frac{\rr_+}{\rr_0}+\frac{\rr_-}{\rr_0}\right)\, , \quad \b=\frac{1}{2}\left(\frac{\rr_+}{\rr_0}-\frac{\rr_-}{\rr_0}\right)\, , \quad \o:=\arg{(-v_2, v_1)}\,,\eeqa
 where we have let, for short,
\beqano\rr_0:=\|{\mathbf v}_0\|\, , \quad \rr_\pm:=\|{\mathbf v}\pm {\mathbf v}_0\|\, . \eeqano

\nl
Regarding $\rr_0$ as a fixed external parameter and calling $p_\l$, $p_\b$, $p_\o$ the generalized momenta associated to $\l$, $\b$ and $\o$, it turns out that
the Hamiltonian~\equ{2Cold}, written in the coordinates $(p_\l, p_\b, \l, \b)$ is independent of $\o$  and has the expression
\beqa{H2C}\JJ(p_\l, p_\b, p_\o, \l, \b, \rr_0)&=&\frac{1}{\l^2-\b^2}\Big[\frac{p^2_\l(\l^2-1)}{2 \rr_0^2}+\frac{p^2_\beta(1-\beta^2)}{2 \rr_0^2}+\frac{p_\o^2}{2 \rr_0^2}\left(\frac{1}{1-\beta^2}\right.\nonumber\\
&&\left.+\frac{1}{\l^2-1}\right)-\frac{(\mm_++\mm_-)\l}{\rr_0^2}+\frac{(\mm_+-\mm_-)\beta}{\rr_0^2}\Big]\, . 
\eeqa
It follows that the solution $W$ of Hamilton--Jacobi equation
\beqa{HJNEW}\JJ(W_\l, W_\b, p_\o, \l, \b, \rr_0)=h\eeqa
can be searched of the form
\beqano W(\l, \b, p_\o, \rr_0, h)=W^\ppu(\l,  p_\o, \rr_0, h)+W^\ppd(\b, p_\o, \rr_0, h)\eeqano
and~\equ{HJNEW} separates completely as
\beqa{split}\cF^\ppu(W^\ppu_\l,\l,{p_\o}, \rr_0, h)+\cF^\ppd(W^\ppd_\b,\b,{p_\o}, \rr_0, h)=0\eeqa
with  $\cF^\ppu$, $\cF^\ppd$  defined via~\equ{H2C}--\equ{HJNEW}.

\nl
The identity~\equ{split} implies that
there must exist a function $\EE$, which we call {\it Euler integral}, depending on  $({p_\o}, \rr_0, h)$ only,  such that
\beqano\cF^\ppu(p_\l,\l,{p_\o}, \rr_0, h)=-\cF^\ppd(p_\b,\b,{p_\o}, \rr_0, h)=\EE({p_\o}, \rr_0, h)\quad \forall\ (p_\l, p_\b, \l, \b)\, . \eeqano 
After elementary computations, one find that, in terms of the initial position--impulse coordinates, the Euler Integral
\beqa{Es}
\EE=\frac{1}{2}\left(\cF^\ppu-\cF^\ppd\right)
\eeqa has the expression in~\equ{G1}, when written in the original coordinates.

\paragraph{The ``asymmetric'' case}\label{The ``asymmetric'' case}
We are interested to find the expression of the Euler integral~\equ{G1} when 2CP is written in the form
\beqa{newH2C}
 \JJ=\frac{\|{{\mathbf y}}\|^2}{2 }
-\frac{1}{\|{\mathbf x}\|}-\frac{\cM'}{\|{\mathbf x}'-{{\mathbf x}}\|}%=:\JJ_0+\ovl\m \JJ_1
\eeqa
namely, when the two centres are in ``asymmetric positions'', $\mathbf 0$, $\mathbf x'$.
As we shall see, in that case we have
  \beqa{ENEW}\EE=\|{\mathbf M}\|^2-{\mathbf x}'\cdot {\mathbf L}+\cM'\frac{({\mathbf x}'-{\mathbf x})\cdot {\mathbf x}'}{\|{\mathbf x}'-{\mathbf x}\|}\eeqa
  where 
\beqa{CL}
{\mathbf M}:={\mathbf x}\times {\mathbf y}\, , \qquad {\mathbf L}:={{\mathbf y}}\times{\mathbf M}-\frac{{{\mathbf x}}}{\|{{\mathbf x}}\|}=\ee {\mathbf P}\eeqa 
are the  {\it angular momentum} and the {\it eccentricity vector}   associated to the Kepler Hamiltonian~\equ{Kep}
with $\ee$ and ${\mathbf P}$ being the eccentricity and the perihelion direction ($\|{\mathbf P}\|=1$).
Notice that $\JJ$ reduces to a Kepler Hamiltonian in two cases: either for ${\mathbf x}'=\mathbf 0$, in which case, as in the symmetric case above, $\EE$ reduces to $\|{\mathbf M}\|^2$, or for $\cM'=0$. The latter case is more interesting to us, as $\JJ$ and $\EE$ become, respectively, 
$\JJ_0$ in~\equ{Kep}
and 
\beqa{EEE}\EE_0&=&\|{\mathbf M}\|^2-{\mathbf x}'\cdot {\mathbf L}\eeqa
with $\EE_0$ being -- as well expected -- a combination of first integrals of $\JJ_0$.
\\
 To prove~\equ{ENEW}--\equ{CL}, %we let
%\beqano\widehat{\JJ}(\widehat{\mathbf y},\widehat{\mathbf x}, \widehat{\mathbf x}'):=\frac{1}{\mm}{\JJ}(\mm \widehat{\mathbf y}, \widehat{\mathbf x}, \widehat{\mathbf x}')=\frac{\|  {\widehat{\mathbf y}}\|^2}{2}-\frac{1}{\|{\widehat{\mathbf x}}\|}-\frac{\cM'}{\|{\widehat{\mathbf x}}-{\widehat{\mathbf x}'}\|}\eeqano
%and then 
we change, canonically,
\beqano{\mathbf x}'=2{\mathbf v}_0\, , \quad {\mathbf x}={\mathbf v}_0+{\mathbf v}\, , \quad {\mathbf y}'=
\frac{1}{2}({\mathbf u}_0-{\mathbf u})
\, , \quad  {\mathbf y}={\mathbf u} \eeqano
(where ${\mathbf y}'$, ${\mathbf u}_0$ denote the generalized impulses conjugated to ${\mathbf x}'$, ${\mathbf v}_0$, respectively)
we reach the Hamiltonian $\JJ$ in~\equ{2Cold}, with $\mm_+=1$, $\mm_-=\cM'$.  Turning back with the transformations, one sees that the function $\EE$ in~\equ{G1} takes the expression
\beqano
{\EE}&:=&\left\|\left({\mathbf x}-\frac{{\mathbf x}'}{2}\right)\times {\mathbf y}\right\|^2+\frac{1}{4}({\mathbf x}'\cdot {\mathbf y})^2+{\mathbf x}'\cdot\left({\mathbf x}-\frac{{\mathbf x}'}{2}\right)\left(\frac{1}{\|{\mathbf x}\|}-\frac{\cM' }{\|{\mathbf x}'-{\mathbf x}\|}\right)\, . 
\eeqano
and we rewrite it as
\beqano
\EE=\EE_0+\EE_1+\EE_2\eeqano
with 
\beqano
&&\EE_0:=\|{\mathbf M}\|^2-{\mathbf x}'\cdot {\mathbf L}\, , \quad 
\EE_1:= \cM'\frac{({\mathbf x}'-{\mathbf x})\cdot {\mathbf x}'}{\|{\mathbf x}'-{\mathbf x}\|}\nonumber\\
&& \EE_2:=\frac{\|{\mathbf x}'\|^2}{2}\left(\frac{\|{{\mathbf y}}\|^2}{2 }
-\frac{1}{\|{\mathbf x}\|}-\frac{\cM'}{\|{\mathbf x}'-{{\mathbf x}}\|}\right)\eeqano
where ${\mathbf M}$, ${\mathbf L}$
are as  in~\equ{CL}.
  Since $\EE_2$ is  itself an integral for ${\JJ}$, we can neglect it and rename 
    \beqa{cal G new}
\EE:=\EE_0+\EE_1\eeqa
the Euler integral to ${\JJ}$. Namely,
\beqa{commutationOLD}
\Big\{\JJ, \EE\Big\}=0\, .
\eeqa

\paragraph{A  set of canonical coordinates which lets  $\JJ$ and $\EE$ in 2 degrees--of--freedom}\label{coordinates}

We describe a set of canonical
coordinates, which we denote as $\cK$, which we shall use for our analysis of the Euler Hamiltonian~\equ{newH2C} and its integral $\EE$~\equ{cal G new}. This set of coordinates puts $\JJ$ and $\EE$ in two degrees--of--freedom (represented by the couples $(\L, \ell)$, $(\GG, {\rm g})$ below), precisely like the classical  ellipsoidal coordinates~\equ{lambdabeta} do, both in the spatial and planar case.

\nl
We consider, in the region of $(\by, \bx)$ 
 where $\JJ_0$ in~\equ{Kep} takes negative values and the ellipse ${\mathbb E}(\by, \bx)$ it generates starting from any initial datum $(\by, \bx)$ in this region is not a circle. Denote as: 
 \begin{itemize}
 \item[{\tiny\textbullet}]
 $a$  the {\it semi--major axis}; 
 % \item[{\tiny\textbullet}]
 %$\ee$  the {\it eccentricity}; 
  \item[{\tiny\textbullet}] ${\mathbf P}$, with $\|{\mathbf P}\|=1$, the direction of perihelion, assuming the ellipse is not a circle;
    \item[{\tiny\textbullet}] $\ell$: the mean anomaly, defined, mod $2\p$, as the area of the elliptic sector spanned by ${\mathbf x}$ from ${\mathbf P}$, normalized to $2\p$.
     \end{itemize}
     Finally,
     \begin{itemize} 
%  \item[{\tiny\textbullet}]  the {\it true anomaly} $\n$, defined as $$\n=\arg(\cos\xi-\ee, \sqrt{1-\ee^2}\sin \xi)$$ with 
 %  \item[{\tiny\textbullet}]  the {\it eccentric anomaly} $\xi$, solving the {\it Kepler equation} $\xi-\ee\sin\xi=\ell$.
  %  \item[{\tiny\textbullet}] $\ee(\L,\GG)=\sqrt{1-\GG^2}$ is the {\it eccentricity};
 %    \item[{\tiny\textbullet}]    the quantity $\varrho=1-\ee\cos\xi$ corresponding to the ratio $\frac{\rr}{a}$.
     \item[{\tiny\textbullet}] given three vectors ${\mathbf u}$, ${\mathbf v}$ and ${\mathbf w}$, with ${\mathbf u}$, ${\mathbf v}\perp{\mathbf w}$, we denote as $\a_{\mathbf w}({\mathbf u}, {\mathbf v})$ the oriented angle  from ${\mathbf u}$ to ${\mathbf v}$ relatively to the positive orientation established by ${\mathbf w}$.
    \end{itemize}

%\nl 
%Next, we introduce the following notations.
%\begin{itemize}
%     \item[{\tiny\textbullet}]  If ${\mathbf i}$, ${\mathbf k}\in \real^3$, with  ${\mathbf i}\perp {\mathbf k}$ and  ${\mathbf i}$, $\mathbf k$ not necessarily of length 1,  by $\FF\sim({\mathbf i},\cdot,{\mathbf k})$, we mean  the orthonormal frame  $\FF=(\frac{{\mathbf i}}{\|{\mathbf i}\|}, \frac{{\mathbf k}\times {\mathbf i}}{\|{\mathbf k}\times {\mathbf i}\|},\frac{{\mathbf k}}{\|{\mathbf k}\|})$.
%     \item[{\tiny\textbullet}] 
%Given a couple $(\FF,\FF')$ of orthonormal frames, where $\FF\sim({\mathbf i},\cdot,{\mathbf k})$, $\FF'\sim({\mathbf i}',\cdot,{\mathbf k}')$, with ${\mathbf k}\times{\mathbf k}'\ne {\mathbf 0}$, we write $$\FF\to^{(\YY, \XX,  x)}\FF'$$ if ${\mathbf i}'={\mathbf k}\times{\mathbf k}'$ and
%$$\YY={\mathbf k}'\cdot \frac{\mathbf k}{\|{\mathbf k}\|}\, , \qquad \XX=\|{\mathbf k}'\|\, , \qquad  x=\a_{\mathbf k}({\mathbf i},{\mathbf i'})$$ 
%where $\a_{\mathbf k}({\mathbf i},{\mathbf i'})$ is the oriented angle  ${\mathbf i}$ to ${\mathbf i}'$, with respect to the counterclockwise orientation established by ${\mathbf k}$. 

%\end{itemize}

\nl
We fix an arbitrary (``inertial'') frame 
\beqano\FF_0:\qquad {\mathbf i}=\left(
\begin{array}{lll}
1\\
0\\
0
\end{array}
\right)\, , \qquad {\mathbf j}=\left(
\begin{array}{lll}
0\\
1\\
0
\end{array}
\right)\, , \qquad {\mathbf k}=\left(
\begin{array}{lll}
0\\
0\\
1
\end{array}
\right)\eeqano in $\real^3$,
and denote as
\beqano{\mathbf M}={\mathbf x}\times {\mathbf y}\, , \quad {\mathbf M}'={\mathbf x}'\times {\mathbf y}'\, , \quad {\mathbf C}={\mathbf M}'+{\mathbf M}\, , \eeqano
where ``$\times$'' denotes skew--product in ${\real}^3$.
Observe the following relations
\beqa{orthogonality}{{\mathbf x}'}\cdot{{\mathbf C}}={{\mathbf x}'}\cdot{\big({\mathbf M}+{\mathbf M}'\big)}={{\mathbf x}'}\cdot{{\mathbf M}}\, , \qquad {\mathbf P}\cdot {\mathbf M}=0\, , \quad \|{\mathbf P}\|=1\, . \eeqa
Assume that the ``nodes''
\beqa{nodes}&&{\mathbf n}_1:={\mathbf k}\times {\mathbf C}\, , \quad {\mathbf n}_2:={\mathbf C}\times {\mathbf x}'\, , \quad {\mathbf n}_3:={\mathbf x}'\times {\mathbf M}\eeqa
do not vanish. 
We   define the coordinates \beqano{\cK}=(\ZZ, \CC, \Theta, \GG, \RR  , \L, \zeta, g, \vartheta, {\rm g}, \rr , \ell)\eeqano
via the following formulae.

\beqa{belline}
\arr{
 \dst {\rm Z}:={\mathbf C}\cdot {\mathbf k}\\ \\
 \dst \CC:=\|{\mathbf C}\|\\ \\
  \dst  \RR:=\frac{{\mathbf y}'\cdot {\mathbf x}'}{\|{\mathbf x}'\|}\\\\
\dst  \L=\sqrt{ a}\\\\
 \dst\GG:=\|{\mathbf M}\|\\\\
 \dst \Theta:=\frac{{\mathbf M}\cdot {\mathbf x}'}{\|{\mathbf x}'\|}\\
 }\qquad\qquad \arr{
 \dst {\rm z}:=\a_{{\mathbf k}}({\mathbf i}, {\mathbf n}_1)\\ \\
 \dst g:=\a_{{\mathbf C}}({\mathbf n}_1, {\mathbf n}_2)\\ \\
  \dst \rr:=\|{\mathbf x}'\|\\ \\
\ell:=\textrm{\rm mean anomaly of } {\mathbf x}\ {\rm on\ \mathbb E}\\ \\
{\rm g}:=\a_{{\mathbf M}}({\mathbf n}_3, {\mathbf M}\times {\mathbf P})\\\\
  \dst\vartheta:=\a_{{\mathbf x}'}({\mathbf n}_2, {\mathbf n}_3)
 }
\eeqa

 % fino qui
    \nl
    The canonical character of ${\cK}$ has been discussed in~\cite{pinzari20b}. In the planar case,  the coordinates~\equ{belline} reduce to the 8 coordinates
 \begin{eqnarray}\label{coord}
\arr{\dst\CC=\|\bx\times \by+\bx'\times \by'\|\\\\
\dst\GG=\|\bx\times \by\|\\\\
\dst\RR=\frac{\mathbf y'\cdot \mathbf x'}{\|\mathbf x'\|}\\\\
\dst\L= \sqrt{ a}
}\qquad\qquad \arr{\dst\g=\a_{\mathbf k}(\mathbf i, \mathbf x')+\frac{\p}{2}\\\\
\dst {\rm g}=\a_{\mathbf k}({\mathbf x'},\mathbf P)+\p\\\\
\dst \rr=\|\mathbf x'\|\\\\
\dst\ell={\rm mean\ anomaly\ of\ {\mathbf x}\ in\ \mathbb E}
}
\end{eqnarray}

\nl
 Using the formulae in the previous section, we provide the expressions of $\JJ$ in~\equ{newH2C} and and $\EE$ in~\equ{EEE} in terms of ${\cK}$:

\beqa{EEJJ}
\JJ(\L, \GG, \Theta, \rr , \ell, {\rm g})&=&-\frac{1}{2\L^2}-\frac{\cM'}{\sqrt{{\rr }^2+2\rr  a\sqrt{1-\frac{\Theta^2}{\GG^2}} {\rm p}+{a}^2\varrho^2}}\nonumber\\
&=:&\JJ_{0}+\JJ_{1}\nonumber\\
\EE(\L, \GG, \Theta, \rr , \ell, {\rm g})&=&\GG^2+\rr \sqrt{1-\frac{\Theta^2}{\GG^2}}\sqrt{1-\frac{\GG^2}{\L^2}}\cos{\rm g}\nonumber\\
 & +&\cM'\rr \frac{{\rr }+a\sqrt{1-\frac{\Theta^2}{\GG^2}} {\rm p}}{\sqrt{{\rr }^2+2\rr  a\sqrt{1-\frac{\Theta^2}{\GG^2}}{\rm p}+{a}^2\varrho^2}}\nonumber\\
&=:&\EE_{0}+\EE_{1}
\eeqa
and, if $\xi=\xi(\L, \GG, \ell)$ is the {\it eccentric anomaly}, defined as the solution of
{\it Kepler equation}
 \beqa{Kepler Equation}\xi-\ee(\L,\GG)\sin\xi=\ell\eeqa
and $a=a(\L)$ the {\it semi--major axis}; $\ee=\ee(\L, \GG)$, the {\it eccentricity} of the ellipse,
$\varrho=\varrho(\L, \GG, \ell)$, ${\rm p}={\rm p}(\L, \GG, \ell, {\rm g})$ are defined as
\beqa{p}
a(\L)&=&{\L^2}\nonumber\\
\ee(\L, \GG)&:=&\sqrt{1-\frac{\GG^2}{\L^2}}\nonumber\\
 \varrho(\L, \GG, \ell)&:=&1-\ee(\L, \GG)\cos\xi(\L, \GG, \ell)\nonumber\\
{\rm p}(\L, \GG, \ell, {\rm g})&:=&(\cos\xi(\L, \GG, \ell)-\ee(\L, \GG))\cos{\rm g}-\frac{\GG}{\L}\sin\xi(\L, \GG, \ell)\sin{\rm g}\, . \eeqa
The angle
\beqa{true anomaly}
\n(\L, \GG, \ell):=\arg\left(\cos\xi(\L, \GG, \ell)-\ee(\L, \GG),\frac{\GG}{\L}\sin\xi(\L, \GG, \ell) \right)
\eeqa
is usually referred to as {\it true anomaly}, so one recognises that ${\rm p}(\L, \GG, \ell, {\rm g})=\varrho \cos(\n+{\rm g})$.

\nl
Observe that $\EE$ and $\JJ$ in~\equ{EEJJ} do not depend on $\CC$, $\ZZ$, $\zeta$, $\g$, $\RR$, $\vartheta$, while the Hamiltonians~\equ{secular},  do not depend on  $\ZZ$, $\zeta$, $\g$,  $\ell$.

\nl
The details on the derivation of the formulae in~\equ{EEJJ} may be found in~\cite{pinzari20b}.

\paragraph{Renormalizable integrability}\label{The secular Euler Hamiltonian}

In this section we review the property of {\it renormalizable integrability} pointed out in~\cite{pinzari19}.

 \nl
 We consider the function $\UU_\b$ in~\equ{Usb} with $\b=1$, which is given by
 \beqa{h1} \UU(\rr, \L, \Theta, \GG, {\rm g})=\frac{1 }{2\p}\int_0^{2\p} 
 \frac{d\ell}{\sqrt{{\rr }^2+2\rr  a\sqrt{1-\frac{\Theta^2}{\GG^2}} {\rm p}+{a}^2\varrho^2}}
% \frac{d\ell}{\|{\mathbf x}'(\rr, \L, \Theta, \GG, {\rm g})-{\mathbf x}(\rr, \L, \Theta, \ell, \GG, {\rm g})\|}
\eeqa 
and the function $$\EE_0=\GG^2+\rr \sqrt{1-\frac{\Theta^2}{\GG^2}}\sqrt{1-\frac{\GG^2}{\L^2}}\cos{\rm g}$$ in~\equ{EEJJ}.
 These two functions have the following remarkable properties:
 \begin{itemize}
 \item[($\cP_1$)]  they have one effective degree--of--freedom, as they depend on one conjugated couple of coordinates: the couple $(\GG, {\rm g}$);
 \item[($\cP_2$)] they Poisson--commute: \beqa{commutation}\Big\{\UU,\ \EE_0\Big\}=0\, . \eeqa
 \end{itemize}
Relation~\equ{commutation} can be proved taking 
the $\ell$--average of~\equ{commutationOLD}, and exploiting that $\JJ_0$ depends only on $\L$; see~\cite{pinzari19}. 
The following definition relies precisely with this situation.

 \begin{definition}[\cite{pinzari19}]\label{def: renorm integr}\rm Let $h$, $g$ be two  functions
of the form
\beqa{HJ}h(p, q, y, x)=\widehat h({\rm I}(p,q), y, x)\, , \qquad g(p, q, y, x)=\widehat g({\rm I}(p,q), y, x)\eeqa
where 
\beqa{D}(p, q, y, x)\in \cD:=\cB\times U\eeqa
with $ U\subset \real^2$, $\cB\subset\real^{2n}$ open and connected, $(p,q)=$ $(p_1$, $\cdots$, $p_n$, $q_1$, $\cdots$, $q_n)$   {conjugate} coordinates with respect to the two--form $\o=dy\wedge dx+\sum_{i=1}^{n}dp_i\wedge dq_i$ and ${\rm I}(p,q)=({\rm I}_1(p,q), \cdots, {\rm I}_n(p,q))$, with
\beqano{\rm I}_i:\ \cB\to \real\, , \qquad i=1,\cdots, n\eeqano
pairwise Poisson commuting:
\beqa{comm}\big\{{\rm I}_i,{\rm I}_j\big\}=0\qquad \forall\ 1\le i<j\le n\qquad i=1,\cdots, n\, . \eeqa
 We say that $h$ is {\it renormalizably integrable via $g$} if there exists a function  \beqano\widetilde h:\qquad {\rm I}(\cB)\times g(U)\to \real\, ,  \eeqano such that
\beqa{renorm}h(p,q,y,x)=\widetilde h({\rm I}(p,q), \widehat g({\rm I}(p,q),y,x))\eeqa
for all $(p, q, y, x)\in \cD$.
\end{definition}
\begin{proposition}[\cite{pinzari19}]\label{rem}
If $h$  is renormalizably integrable via $g$, then:\item[{\rm (i)}] 
${\rm I}_1$, $\cdots$, ${\rm I}_n$  are first integrals to $h$ and $g$;
\item[{\rm (ii)}] 
$h$ and $g$ Poisson commute. \end{proposition}

\begin{proposition}[\cite{pinzari19}]\label{main prop}
$\UU$ is renormalizably integrable via ${\rm E}_0$. Namely,
there exists a function $\FF$ such that
\beqano\UU({\rm r}, \L, \Theta, {\rm G}, {\rm g})=\FF\big({\rm r}, \L, \Theta, {\rm E}_0({\rm r}, \L, \Theta, {\rm G}, {\rm g})\big)\, . \eeqano

\end{proposition}

%We  recall a result from {referenza}.
%Let
%\beqa{C}\cC:\qquad(\L,\ell,u,v)\in \cA\times\torus\times V\to (\underline{\mathbf y}, \underline{\mathbf x})=({\mathbf y}', {\mathbf y}, {\mathbf x}', {\mathbf x})\in (\real^3)^4\eeqa
% where $\cA$ is domain\footnote{{By ``domain'' we mean an open and connected set in ${\mathbb K}=\real^m, \complex^m$.}} in $\real$,   $V$ is a domain in $\real^{10}$, $\torus:=\real/(2\p\integer)$, $(u, v)=\big((u_1, u_2, u_3, u_4, u_5)$, $(v_1, v_2, v_3, v_4, v_5)\big)$, 
%be a change of coordinates  which verifies
%$$d{\mathbf y}'\wedge d{\mathbf x}'+d{\mathbf y}\wedge d{\mathbf x}=d\L\wedge d\ell+d u\wedge d v$$
%and
% \beq{2B}
%\left(\frac{\|{\mathbf y}\|^2}{2}-\frac{1}{\|{\mathbf x}\|}\right)\circ\cC=
%-\frac{^3\cM^2}{2}=:{\rm h}_{\rm Kep}(\L)\, , \eeq 
%where $$, $\cM$ are fixed having assimed that the image of $\cC$ in~\equ{C}
%is a domain of  $(\mathbf y, \mathbf x)$ where the left hand side of~\equ{2B} takes negative values.  Assume also that the  ellipse ${\mathbb E}$ generated by the two--body Hamiltonian~\equ{2B} has non--vanishing eccentricity. 
%We call any map of such kind {\it partial Kepler map}. 

\nl
%We  consider the ellipse generated by the ``Kepler Hamiltonian'' at left hand side in~\equ{2B} and denote as  ${\rm e}:=\sqrt{1-\frac{{\rm G}^2}{\L^2}}$ its eccentricity, where ${\rm G}:=\|{\mathbf M}\|$; ${\mathbf P}$, with $\|{\mathbf P}\|=1$ the direction of perihelion, assuming ${\rm e}\ne 0$. Recall that the eccentricity vector ${\mathbf L}$ in~\equ{CL} is related to $\ee$ and ${\mathbf P}$ via  ${\mathbf L}=\ee{\mathbf P}$, so that  the function $\EE_0$ in~\equ{Kep} becomes
%\beqa{cG}{{\rm E}_0}:={\rm G}^2-{\rm m}^2{\rm M} {\rm e}\,{\mathbf x}'\cdot {\mathbf P}\eeqa

\nl
The proof of Proposition~\ref{main prop} is based on $\cP_1\div\cP_2$ above. Below, we list some consequences.
\begin{itemize}
\item[(i)] If $\FF_{\EE_0}\ne 0$, the  time laws of $({\rm G}, {\rm g})$ under $\UU$ or $\EE_0$ are basically (i.e., up to a change of time) the same;

\item[(ii)] Motions of $\EE_0$ corresponding to level sets for which  $\FF_{\EE_0}= 0$ are fixed points curves to $\UU$ (``frozen orbits''). In~\cite{pinzari19} we provided an example of frozen orbit of $\UU$ in the spatial case, for $\rr\ll 1$;

\item[(iii)] $\UU$ and $\EE_0$ have the same action--angle coordinates;

\item[(iv)] $\FF$ may have several expressions, as well as $\UU$, which is defined via a quadrature. Two different representation formulae have been proposed in~\cite{pinzari19} and~\cite{pinzari20a}. 

\end{itemize}

\nl
In the next section, we investigate the dynamical properties of $\EE_0$ for the planar case ($\Theta=0$).

%\nl
% We consider the Taylor and  Taylor--Fourier expansions
%\beqano \UU(\rr, \L, \Theta, \GG, {\rm g})=\sum_{n=0}^{+\infty} \rr^{n} \UU_{n}(\L, \Theta, \GG, {\rm g})=\sum_{n=0}^{+\infty}\sum_{m=0}^\infty \rr^{n} \UU_{nm}(\L, \Theta, \GG)\cos m{\rm g}\, . \eeqano
%In {referenza} we have discussed the following consequence of Theorem~\ref{partial integral}. 
% \begin{proposition}\label{HarrK}
% \beqa{Hprop1} \UU_{nm}(\L, \Theta, \GG)\equiv 0\qquad {\rm if}\qquad m\ge n-1\, , \quad \forall \ n\ge 1\, . \eeqa
% In particular, the term $\UU_{1}$ vanishes identically and $\UU_{2}$, called {the} dipolar term,  Poisson--commutes with  $\GG$ (``Harrington property'').
% \end{proposition}

 \paragraph{The phase portrait of $\EE_0$ in the planar case}\label{Dynamical consequences}
Here we fix $\L=1$, $\Theta=0$.
 For  $\rr\in (0,2)$, the function $\EE_0({\rm g}, \GG)$ has a minimum, a saddle and a maximum, respectively at
\beqano{\mathbf P}_{-}=(\pm\p,0)\, , \qquad {\mathbf P}_0=(0,0)\, , \qquad  {\mathbf P}_+=\left(0,\sqrt{1-\frac{\rr^2}{4}}\right)\eeqano
where it takes the values, respectively,
\beqano\cE_-=-\rr\, , \qquad\cE_0=\rr\, , \quad \cE_+=1+\frac{\rr^2}{4}\, . \eeqano
Thus,  the level sets in~\equ{level sets} are non--empty only for \beqa{cal J}\cE\in %[\cE_-,\cE_+]=
\left[-\rr, 1+\frac{\rr^2}{4}\right]\, . \eeqa  
We
denote as $\cS_0$,  the level set through  the saddle ${\mathbf P}_0$. When $\GG=1$, $\EE_0$ takes the value $1$ for all ${\rm g}$ and we denote as $\cS_1$ the level curve with $\cE=1$. The equations of $\cS_0$, $\cS_1$ are, respectively:
\beqa{separatrices}
&&\cS_0(\rr)=\Big\{
({\rm g}, \GG):\ \GG^2+\rr \sqrt{1-\GG^2}\cos{\rm g}=\rr
\Big\},
\nonumber\\
&&\cS_1(\rr)=\Big\{\GG=\pm 1\Big\}\cup\Big\{
\GG=\pm\sqrt{1-\rr^2\cos^2{\rm g}}
\Big\}.
\eeqa
%Observe that 
%$\cS_0$ has the geometrical meaning of ``collisional manifold, since, solving for $\rr$ and then going back with the rescalings~\equ{delta}, one has
%\beqano\rr=\frac{a(1-\ee^2)}{1-\ee\cos{\rm g}}\eeqano
%which says that the point, on the Keplerian ellipse with semi--axis $a$ and eccentricity $\ee$ corresponding to the true anomaly $ {\rm g}+\p$ has radial distance from the sun equal to $\rr$.
$\cS_1$ is composed of two branches,
which will be referred to as ``horizontal'', ``vertical'', respectively, 
transversally intersecting at
$(\pm \frac{\p}{2}, 1)$, with ${\rm g}$ mod $2\p$.
Note that, when $0< \rr< 1$, the vertical branch is defined for all ${\rm g}\in\torus$; when $\rr>1$,  its domain in ${\rm g}$ is made of two disjoint neighbourhoods of $\pm \frac{\p}{2}$. \\
When $\rr>2$, the saddle ${\mathbf P}_0$ and its manifold $\cS_0$ do not exist, ${\mathbf P}_-=(\p, 0)$ is still a minimum, while the maximum becomes ${\mathbf P}_+=(0, 0)$. The manifold $\cS_1$ still exists, with the vertical branch closer and closer, as $\rr\to+\infty$, to the portion of straight ${\rm g}=\pm\frac{\p}{2}$ in the strip $-1\le \GG\le 1$.
In this case the admissible values for $\cE$ are
\beqano\cE\in %[\cE_-,\cE_+]=
\left[-\rr, \rr\right]\, . \eeqano
It is worth mentioning~\cite{pinzari20b} that, when $0<\rr<2$, the motions generated by $\EE_0$ along $\cS_0(\rr)$ can be explicitly computed, and are given by
\beqano
\arr{
\GG(t)=\frac{\s\L}{\cosh\s\L(t-t_0)},\\\\
\dst{\rm g}(t)=\pm\cos^{-1}\frac{1-\frac{\a^2}{\cosh^2\s\L(t-t_0)}}{\sqrt{1-\frac{\s^2}{\cosh^2\s\L(t-t_0)}}}\,,  
}\eeqano
where
\beq{conditions}\s^2:=\rr(2-\rr)\, , \quad \a^2:=2-\rr\, , \quad \rr\in (0,2)\, . \eeq
These motions -- which have a remarkable similitude with the separatrix motions of the classical pendulum -- are however meaningless for $\UU$, which is singular on $\cS_0(\rr)$.

\nl
The scenario is depicted in Figure~\ref{fig:fig1} {to which we refer for further qualitative details}.

\section{Perihelion librations in the three--body problem}\label{Asymptotic  action--angle coordinates}

{In this section we review the results of~\cite{pinzari19},~\cite{pinzari20a},~\cite{pinzari20b}.}

\noindent
As mentioned in the introduction, the proof of Theorem~\ref{main} is based on two ingredients: a normal form theory well designed around the Hamiltonians~\equ{secular}, where no non--resonance condition is required, and a set of action--angle--like coordinates which  approximate well the natural action--angle coordinates of $\EE_0$ when $\rr$ is large. In this section we briefly summarise the procedure. Full details may be found in~\cite{pinzari20a}.

\paragraph{A normal form theory without small divisors}

{W}e describe a procedure for eliminating the angles\footnote{Note  that the procedure described in this section does not seem to be related to~\cite[\S 6.4.4]{arnoldKS06}, for the lack of slow--fast couples.} $\bm\f$ at high orders, given Hamiltonian of the form
 \beqa{h0f0}\HH({\mathbf I}, \bm{\f}, {\mathbf p}, {\mathbf q}, y, x)=\hh({\mathbf I},{\mathbf J}({\mathbf p}, {\mathbf q}), y)+f({\mathbf I}, \bm{\f}, {\mathbf p}, {\mathbf q}, y, x)\eeqa
%$$\hh(y,{\mathbf I},{\mathbf J})={\mathbf J}({\mathbf I}, y)+{\mathbf P}(y,{\mathbf I},{\mathbf J})\, , \qquad {\mathbf P}(y,{\mathbf I},0)\equiv0$$
which we assume to be holomorphic on the neighbourhood
\beqano{\mathbb P}_{\r, s, \d, r, \xi}={\mathbb I}_\r  \times {{\mathbb T}}^n_s\times {\mathbb B}_{\d}\times{\mathbb Y}_r\times  {\mathbb X}_\xi\supset {\mathbb P}={\mathbb I}  \times {{\mathbb T}}^n\times {\mathbb B}\times{\mathbb Y}\times  {\mathbb X}\, , \eeqano
for suitable $\r$, $s$, $\d$, $r$, $\xi>0$ and
\beqano{\mathbf J}({\mathbf p}, {\mathbf q})=(p_1q_1,\cdots, p_mq_m)\, . \eeqano
Here, ${\mathbb I}\subset {\mathbb R}^n$, ${\mathbb B}\subset {\mathbb R}^{2m}$, ${\mathbb Y}\subset {\mathbb R}$, ${\mathbb X}\subset {\mathbb R}$ are open and connected; $\mathbb T=\mathbb R/(2\p\mathbb Z)$ is the standard torus, and we have used the common notation $A_r:=\bigcup_{x\in A}B_r(x)$, where $B_r(x)$ is the complex open ball centered in $x$ with radius $r$.

\nl
We denote as ${\cal O}_{\r, s, \d, r, \xi}$ the set of complex holomorphic functions \beqano\phi:\quad {\mathbb P}_{\hat\r, \hat s, \hat \d, \hat r, \hat\xi}\to {\complex}\eeqano for some  $\hat\r>\r$, $\hat s>s$, $\hat\d>\d$, $\hat r>r$, $\hat\xi>\xi$, equipped with 
 the norm
\beqano\|\phi\|_{\r, s, \d, r, \xi}:=\sum_{k,h,j}\|\phi_{khj}\|_{\r, r,\xi}e^{s|k|}\d^{h+j}\eeqano
where $\phi_{khj}({\mathbf I}, y,  x)$ are the coefficients of the Taylor--Fourier expansion\footnote{We denote as ${\mathbf x}^h:=x_1^{h_1}\cdots x_n^{h_n}$, where ${\mathbf x}=(x_1, \cdots, x_n)\in {\real}^n$ and $h=(h_1, \cdots, h_n)\in {\natural}^n$.}
\beqano\phi=\sum_{k,h,j}\phi_{khj}({\mathbf I},y, x)e^{\ii k s}{\mathbf p}^h {\mathbf q}^j\, , \quad \|\phi\|_{\r, r,\xi}:=\sup_{{\mathbb I}_\r\times{\mathbb Y}_r\times  {\mathbb X}_\xi}|\phi({\mathbf I}, y, x)|\, . \eeqano
If $\phi$ is independent of $x$, we simply write
$\|\phi\|_{\r, r}$ for $\|\phi\|_{\r, r,\xi}$.
If $\phi\in {\cal O}_{\r, s, \d, r, \xi}$, we define its ``off--average'' $\widetilde\phi$ and   ``average'' $\ovl\phi$ as
\beqano&&\widetilde\phi:=\sum_{k,h,j:\atop (k,h-j)\ne (0,0)}\phi_{khj}({\mathbf I}, y,  x)e^{\ii k s}{\mathbf p}^h {\mathbf q}^j\nonumber\\
&&\ovl\phi:=\phi-\widetilde\phi=\frac{1}{(2\p)^n}\int_{[0,2\p]^n}\P_{{\mathbf p}\mathbf q}\phi({\mathbf I}, \bm{\f}, {\mathbf J}({\mathbf p}, {\mathbf q}), y, x)d\bm{\f}
\, , \eeqano
with
\beqano\P_{{\mathbf p}\mathbf q}\phi({\mathbf I}, \bm{\f}, {\mathbf J}({\mathbf p}, {\mathbf q}), y, x):=\sum_{k,h}\phi_{khh}({\mathbf I}, y,  x)e^{\ii k s}{\mathbf p}^h {\mathbf q}^h\eeqano
We decompose
\beqano{\cal O}_{\r, s, \d, r, \xi}={\cal Z}_{\r, s, \d, r, \xi}\oplus {\cal N}_{\r, s, \d, r, \xi}\, . \eeqano where ${\cal Z}_{\r, s, \d, r, \xi}$, ${\cal N}_{\r, s, \d, r, \xi}$
are the ``zero--average'' and the  ``normal'' classes
\beqa{zero average}&&{\cal Z}_{\r, s, \d, r, \xi}:=\{\phi\in {\cal O}_{\r, s, \d, r, \xi}:\quad \phi=\widetilde\phi\}=\{\phi\in {\cal O}_{\r, s, \d, r, \xi}:\quad \ovl\phi=0\}\\
\label{phi independent}&&{\cal N}_{\r, s, \d, r, \xi}:=\{\phi\in {\cal O}_{\r, s, \d, r, \xi}:\quad\phi=\ovl\phi\}=\{\phi\in {\cal O}_{\r, s, \d, r, \xi}:\quad \widetilde\phi=0\}\, . \eeqa
respectively. 
We finally let
$\omega_{y,{\mathbf I},{\mathbf J}}:=\partial_{y,{\mathbf I},{\mathbf J}} \hh$.

\nl
In the following result, no non--resonance condition is required on the frequencies $\o_{\mathbf I}$, which, as a matter of fact, might also be zero.

\begin{theorem}[\cite{pinzari20a}]\label{NFL}
For any $n$, $m$, there exists a number ${\rm c}_{n,m}\ge 1$ such that, for any $N\in {\mathbb N}$ such that the following inequalities are satisfied
\beqa{normal form assumptions}
&&4N{\cal X}\left\|\Im\frac{\o_{\mathbf I}}{\o_y}\right\|_{\r, r}<
s
\, , \quad
4N{\cal X}\left\|\frac{\o_{\mathbf J}}{\o_y}\right\|_{\r, r}<
1\nonumber\\
&&   {\widetilde{\rm c}_{n,m}N\frac{{\cal X}}{\dd}   % \|\frac{1}{\o_y}\|_{\r, r}
 \left\|{f}\right\|_{\r, s, \d, r, \xi}\left\|\frac{1}{\o_y}\right\|_{\r, s, \d, r, \xi} <1}
  %&& {\rm c}^{2}_{n,m}N^{2}\frac{{\cal X}^{2}}{{\rm d}^{2}}  { \left\|\frac{1}{\o_y}\right\|^2_{\r, r} 
%\left\|% \frac
%{f}%{\o_y}
%\right\|_{\r, s, \d, r, \xi}}\left\|% \frac
%\widetilde{f}%{\o_y}
%\right\|_{\r, s, \d, r, \xi}<1 
\eeqa
with ${\rm d}:=\min\big\{\r s, r\xi, {\d}^2\big\}$, ${\cal X}:=\sup\big\{|x|:\ x\in {\mathbb X}_\xi\big\}$, one can find an operator \beqa{PsiN}\Psi_*:\quad {\cal O}_{\r, s, \d, r, \xi}\to \cO_{1/3 (\r, s, \d, r, \xi)}\eeqa
which carries $\HH$ to
\beqano\HH_*=\hh+g_*+f_*\eeqano
where
$g_*\in \cN_{1/3 (\r, s, \d, r, \xi)}$, $f_*\in \cO_{1/3 (\r, s, \d, r, \xi)}$ and, moreover, the following inequalities hold 
\beqa{thesis}
&&\|g_*-\ovl f\|_{1/3 (\r, s, \d, r, \xi)}\le 162\widetilde{\rm c}_{n, m} \frac{{\cal X}}{\rm d}\left\|%\frac{1}{\o_y}\|_{\r, r} \|
{\frac{\widetilde f}{\o_y}}\right\|_{\r, s, \d, r, \xi}\| f\|_{\r, s, \d, r, \xi}\nonumber\\
&&    \|f_*\|_{1/3 (\r, s, \d, r, \xi)}\le \frac{1}{2^{N+1}} \|f\|_{\r, s, \d, r, \xi}\, . \eeqa
The transformation $\Psi_*$ can be obtained as a composition of time--one Hamiltonian flows, and satisfies the following.  If
\beqano({\mathbf I}, \bm{\f}, {\mathbf p}, {\mathbf q}, y, x):=\Psi_*({\mathbf I}_*, \bm{\f}_*, {\mathbf p}_*, {\mathbf q}_*, \RR_*, \rr_*)\eeqano the following uniform bounds hold:
\beqa{phi close to id***}
&&{\rm d}\max\Big\{\frac{|{\mathbf I}-{\mathbf I}_*|}{\r},\ \frac{|\bm{\f} -\bm{\f} _*|}{s},\ 
\frac{|{\mathbf p}-{\mathbf p}_*|}{\d},\ \frac{ |{\mathbf q}-{\mathbf q}_*|}{\d},\ \frac{|y-y_*|}{r},\frac{ |x-x_*|}{\xi} \Big\}\nonumber\\
&&\le\max\Big\{s|{\mathbf I}-{\mathbf I}_*|,\ \r|\bm{\f} -\bm{\f} _*|,\ 
\d |{\mathbf p}-{\mathbf p}_*|,\ \d |{\mathbf q}-{\mathbf q}_*|,\ \xi |y-y_*|,\ r |x-x_*| \Big\}\nonumber\\
&&\leq%\yellow{4\, {\cal X}   % \|\frac{1}{\o_y}\|_{\r, r}
% \| \frac{f}{\o_y}\|_{\r, s, \d, r, \xi}}
 {{19\, {\cal X}   % \|\frac{1}{\o_y}\|_{\r, r}
\left \| \frac{ f}{\o_y}\right\|_{\r, s, \d, r, \xi}
}}%+
 %4\, \frac{{\cal X}^2 N}{\dd\|\o_y\|_{\r, s, \d, r, \xi}}  % \|\frac{1}{\o_y}\|_{\r, r}
 %\left\| \frac{\widetilde f}{\o_y}\right\|_{\r, s, \d, r, \xi}
%{\|f\|_{\r, s, \d, r, \xi}}
 %}
 \, . 
\eeqa
\end{theorem}
Hints on the proof of Theorem~\ref{NFL} may be found in Appendix~\ref{Outline of the proof}.
\paragraph{Asymptotic  action--angle coordinates}

The explicit construction of the action--angle coordinates for $\EE_0$ for any value of $\rr$  and $\Theta$ exhibits elliptic integrals. This is true even in the case $\Theta=0$, in which the phase portrait is, as discussed, explicit. 
As we are interested to the case that $\rr$ is large, we adopt the ``approximate'' solution  of integrating only the leading part of $\EE$. Namely, we replace 
equation~\equ{level sets} with

\beqa{tildecE} \sqrt{1-\GG^2}\cos{\rm g}=\widetilde\cE\,.\eeqa
  We show that, for this case, the action--angle coordinates, denoted as $(\cG, \g)$, are given by
  \beqa{goodtransf} \cG=\widetilde\cE
\, , \quad \g=\t\eeqa
where $\t$ is the time the flows employs to reach the value $(\GG, {\rm g})$ on the level set $\widetilde\cE$, starting from $(\sqrt{1-\widetilde\cE^2}, 0)$ ($(\sqrt{1-\widetilde\cE^2}, \p)$).
Namely, for the Hamiltonian~\equ{tildecE}, the action--angle coordinates  coincide with the energy--time coordinates.
 Indeed, by~\equ{tildecE}, the action variable can be taken to be

 \beqano
 \cG(\widetilde\cE)&=&\arr{
\dst-1+\frac{1}{\p}\int_{-\arccos|\widetilde\cE|}^{\arccos|\widetilde\cE|}\sqrt{1-\frac{\widetilde\cE^2}{\cos^2{\rm g}}}d{\rm g}\quad -1< \widetilde\cE< 0\\\\
\dst1-\frac{1}{\p}\int_{-\arccos\widetilde\cE}^{\arccos\widetilde\cE}\sqrt{1-\frac{\widetilde\cE^2}{\cos^2{\rm g}}}d{\rm g}\quad 0< \widetilde\cE< 1\,.
}
\eeqano
We have defined  $\cG(\widetilde\cE)$ so that $\cG(0)=0$. Then 
the period of the orbit is given  by
 \beqano\cT(\widetilde\cE)=2\p \cG_{\widetilde\cE}(\widetilde\cE)\,.\eeqano
With the change of variable \beqa{w}w=\frac{|\widetilde\cE|}{\sqrt{1-\widetilde\cE^2}}\tan{\rm g}\,,\eeqa we obtain
\beqano
\cT(\widetilde\cE)&=&4|\widetilde\cE|\int_{0}^{\arccos|\widetilde\cE|}\frac{1}{\cos^2{\rm g}}\frac{d{\rm g}}{\sqrt{1-\frac{\widetilde\cE^2}{\cos^2{\rm g}}}}=4\int_{0}^{1}\frac{dw}{\sqrt{1-w^2}}=2\p
\eeqano
which implies~\equ{goodtransf}. \\ 
Looking at the (multi--valued) generating function
\beqano S(\cG, {\rm g})=\int_{P_0(\cG)}^{P_{\rm g}(\cG)}\sqrt{1-\frac{\cG^2}{\cos^2{\rm g}'}}d{\rm g}'\eeqano
(where, as it is standard to do~\cite{arnold89}, 
the integral is computed along the  $\cG$th level set, from $P_0(\cG):=(-\arccos\cG, 0)$ to a prefixed point $P_{\rm g}(\cG)=({\rm g}, \cdot)$ of the level set, so as to make $S(\cG, \cdot)$ continuous) we obtain the  transformation of coordinates
\beqa{Gg} 
 \arr{
\dst\GG=\sqrt{1-\cG^2}\cos\g\\\\
\dst {\rm g}=-\tan^{-1}\left(\frac{1}{\cG}\sqrt{1-{\cG^2}}\sin\g\right)+k\p\\\\ {\rm with}\quad k=\arr{0\ \ {\rm if}\ 0<\cG<1\\
1\ \ {\rm if}\ -1<\cG<0\,.
}%\ \sign\cos{\rm g}=\sign\cG
}
\eeqa
Then, using the coordinates $(\cG, \g)$,  one obtains the expression
\beqano
\EE_0={\cG}\,\rr+(1-\cG^2)\cos^2\g
\eeqano which will be used in the next section.
%\footnote{
%Here we have used the expansion
%\beqano\UU= -\frac{\cM'}{\rr}-\cM'\frac{3}{2}\frac{a}{\rr^2}\sqrt{1-\GG^2}\cos{\rm g}+\OO(\rr^{-3})\, . \eeqano
%Observe the proportionality between $\EE_{0, 1}$ and the term of order $\rr^{-2}$ in the expansion of $\UU$. Such proportionality persists also in the general case: 	\beqano\EE_{0, 1}= \rr \sqrt{1-\GG^2} \sqrt{1-\frac{\Theta^2}{\GG^2}}\cos{\rm g}\eeqano
%and
%\beqano\UU= -\frac{\cM'}{\rr}-\cM'\frac{3}{2}\frac{a}{\rr^2}\sqrt{1-\GG^2}\sqrt{1-\frac{\Theta^2}{\GG^2}}\cos{\rm g}+\OO(\rr^{-3})\, . \eeqano
%Relation~\equ{commutation} actually implies a weaker assertion: that the lower order terms in the respective expansions of $\EE_0$ and $\UU+\frac{\cM'}{\rr}$ do commute. 
% }  of $\JJ_{\rm s}$ with respect to $\rr$, which is
%\beqa{s, trunc}\JJ_{\rm s}(\rr, \L,  \GG,  {\rm g})=-\frac{^3\cM^2}{2}-\frac{\cM'}{\rr}-\cM'\frac{3}{2}\frac{a}{\rr^2}\sqrt{1-\GG^2}\cos{\rm g}+\OO(\rr^{-3})\, , \eeqa

\nl
The proof of Theorem~\ref{main} is a direct application of Proposition~\ref{NFL}. As such, a careful evaluation of the involved quantities is needed. Those evaluations are completely explicit in~\cite{pinzari20a}. For the purpose of this review,  we report below the main ideas while we  skip most of computational details. The reader who is interested in them might consult~\cite{pinzari20a}.

\paragraph{Sketch of proof of Theorem~\ref{main} }
For definiteness, we sketch  the proof of Theorem~\ref{main} for $(0, 0)$. The proof for $(0, \p)$ is similar.
For the purposes of this proof, we let $\widehat\HH_1:=\widehat\HH_\JJ$ and $\widehat\HH_2:=\widehat\HH_0$, where $\widehat\HH_\JJ$ and $\widehat\HH_\JJ$ are as in~\equ{secular}.
It is convenient to rewrite the functions $\widehat\HH_i$  as
\beqano
\widehat\HH_1(\RR, \GG, \rr, {\rm g})&=&\left(\frac{\RR^2}{2}-\frac{1}{\rr}\right)+\frac{\GG^2}{2\rr^2}-\frac{\ovl\b}{\b+\ovl\b}\frac{1}{\rr}\left(\widehat\FF_{\b\varepsilon(\rr)}\Big(\widehat\EE_{\b\varepsilon(\rr)}(\GG, {\rm g})\Big)-1\right)\nonumber\\\nonumber\\
&&-\frac{\b}{\b+\ovl\b}\frac{1}{\rr}\left(\widehat\FF_{-\ovl\b\varepsilon(\rr)}\Big(\widehat\EE_{-\ovl\b\varepsilon(\rr)}(\GG, {\rm g})\Big)-1\right)\nonumber\\
\widehat\HH_2(\RR, \GG, \rr, {\rm g})&=&\left(\frac{\RR^2}{2}-\frac{1}{\rr}\right)+\frac{\GG^2}{2\rr^2}-\frac{\ovl\b}{\b+\ovl\b}\frac{1}{\rr}\left(\widehat\FF_{(\b+\ovl\b)\varepsilon(\rr)}\Big(\widehat\EE_{(\b+\ovl\b)\varepsilon(\rr)}(\GG, {\rm g})\Big)-1\right)
\eeqano
where
\beqano
\varepsilon(\rr):=\su\rr\,,\quad \widehat\EE_\varepsilon(\GG, {\rm g}):=\varepsilon\EE_0(\varepsilon^{-1}, \GG, {\rm g})\, , \quad \widehat\FF_\varepsilon(t):=\varepsilon^{-1}\,\FF(\varepsilon^{-1},\varepsilon^{-1}t)\, .
\eeqano
We next change coordinates via the canonical changes
\beqano\cC_1:\ (\cG, \g)\to (\GG, {\rm g})\, , \quad \cC_2:\ (y, x)\to (\RR, \rr)\eeqano where 
$\cC_1$ is defined as in~\equ{Gg}, with $k=0$, while $\cC_2$ is 
 \beqa{Rr}\arr{\dst\RR(y, x)=\frac{1}{y}\sqrt{ \frac{\cos\x'(x)+1}{1-\cos\xi'(x)}}
 \\\\
 \dst \rr(y, x)=y^2(1-\cos\xi'(x))}\eeqa
where $\xi'(x)$ solves
\beqa{Kepler1}\xi'-\sin\xi'=  x\, . \eeqa
$\cC_2$ has been chosen so that
\beqano \left(\frac{\RR^2}{2}-\frac{1}{\rr}\right)\circ\cC_2=-\frac{1}{2y^2}\, . \eeqano Using the new coordinates, we have
\beqano
\widehat\HH_1&=&-\frac{1}{2y^2}+\frac{1}{\rr(y, x)}\left(\varepsilon(y, x)\frac{(1-\cG^2)}{2}\cos^2\g-\frac{\ovl\b}{\b+\ovl\b}\left(\widehat\FF%^+
_{\b\varepsilon(y, x)}\Big(\widehat\EE_{\b\varepsilon(y, x)}(\cG, \g)\Big)-1\right)\right.\nonumber\\
&&\left.-\frac{\b}{\b+\ovl\b}\left(\widehat\FF_{-\ovl\b\varepsilon(y, x)}\Big(\widehat\EE_{-\ovl\b\varepsilon(y, x)}(\cG, \g)\Big)-1\right)\right)\nonumber\\\nonumber\\
\widehat\HH_2&=&-\frac{1}{2y^2}+\frac{1}{\rr(y, x)}\left(\varepsilon(y, x)\frac{(1-\cG^2)}{2}\cos^2\g\right.\nonumber\\
&&\left.
-
\frac{\ovl\b}{\b+\ovl\b}\frac{1}{\rr(y, x)}\left(\widehat\FF%^+
_{(\b+\ovl\b)\varepsilon(y, x)}\Big(\widehat\EE_{(\b+\ovl\b)\varepsilon(y, x)}(\cG, \g)\Big)-1\right)\right)
\eeqano
having abusively denoted as $\varepsilon(y, x):=\varepsilon\big(\rr(y, x)\big)$ and
\beqa{hatE}\widehat\EE_\varepsilon(\cG, \g):=\cG+\varepsilon\left(1-{\cG^2}\right)\cos^2\g\, . \eeqa
A domain where we shall check holomorphy for $\widehat\HH_i$  is chosen as 
\beqa{domain}{\mathbb D}_{\d, s_0, \sqrt{\a_-},\sqrt{\varepsilon_0}}:={\mathbb Y}_{\sqrt{\a_-}}\times {\mathbb X}_{\sqrt{\varepsilon_0}}\times {\mathbb G}_\d\times {\mathbb T}_{s_0}\eeqa
where
\beqa{yx}
&&{\mathbb Y}:=\Big\{y\in \real:\ 2\sqrt{\a_-}<y<\sqrt{\a_+}\Big\}\, , \quad {\mathbb X}:=\Big\{x\in \real:\ |x-\p|\le \p-2\sqrt{\varepsilon_0}\Big\}\nonumber\\
&&{\mathbb G}:=\Big\{\cG\in \real:\ 1-\d<\cG<1\Big\}\eeqa
where $0<\a_-<\frac{\a_+}{4}$, $0<\varepsilon_0<\frac{\p^2}{4}$, $0<\d<1$. If 
$c_0>0$ is such that   for any $0<\varepsilon_0<1$ and for any  $x\in{\mathbb X}_{\sqrt{\varepsilon_0}}$,
Equation~\equ{Kepler1} has a unique solution $\xi'(x)$ which depends analytically on $x$ and and verifies \beqa{xiineq} |1-\cos \xi'(x)|\ge c_0\varepsilon_0\eeqa
(the existence of such a number $c_0$  is well known)
and if the following inequalities
 are satisfied

\beqa{ass1}&&0<\d\le \frac{1}{4}\, , \qquad C^*(s_0)\d<1\qquad C^*(s_0):=16\left(\sup_{ \torus_{s_0}}|\sin\g|\right)^2\nonumber\\
&& \a_-\varepsilon_0>\frac{4\b^* }{c_0}\eeqa
with $\b^*$ as in~\equ{b*}, then $\widehat\HH_i$ are holomorphic in the domain~\equ{domain}. The proof is based on the explicit evaluation of the function $\widehat\FF_\varepsilon(t)$ for complex values of its arguments, the accomplishment of which is obtained  using the explicit expression of $\widehat\FF_\varepsilon(t)$: see~\cite[Proposition 3.1 and Proposition 3.3]{pinzari20a}.

\nl
We  aim to apply Theorem~\ref{NFL}, with ${\mathbf I}=\cG$, $\bm\f=\g$,
$(y, x)$ as in~\equ{Rr},
$\hh(y)=-\frac{1}{2y^2}$
and, finally

\beqa{f0}
f(\cG, \g, y, x)=\arr{\frac{1}{\rr(y, x)}\left(\varepsilon(y, x)\frac{(1-\cG^2)}{2}\cos^2\g-\frac{\ovl\b}{\b+\ovl\b}\left(\widehat\FF%^+
_{\b\varepsilon(y, x)}\Big(\widehat\EE_{\b\varepsilon(y, x)}(\cG, \g)\Big)-1\right)\right.\\
\left.-\frac{\b}{\b+\ovl\b}\left(\widehat\FF_{-\ovl\b\varepsilon(y, x)}\Big(\widehat\EE_{-\ovl\b\varepsilon(y, x)}(\cG, \g)\Big)-1\right)\right)\quad i=1\\\\
\frac{1}{\rr(y, x)}\left(\varepsilon(y, x)\frac{(1-\cG^2)}{2}\cos^2\g\right.\\
\left.-\frac{\ovl\b}{\b+\ovl\b}\frac{1}{\rr(y, x)}\left(\widehat\FF%^+
_{(\b+\ovl\b)\varepsilon(y, x)}\Big(\widehat\EE_{(\b+\ovl\b)\varepsilon(y, x)}(\cG, \g)\Big)-1\right)\right)\\
\quad i=2
}
\eeqa
As $\hh$ does not depend on $\cG$ and the coordinates ${\mathbf p}$, ${\mathbf q}$ do not exist, in order to apply Theorem~\ref{NFL}, only 
the last condition in~\equ{normal form assumptions} needs to be verified. Direct computations show that such condition is verified provided that $N=[N_0]-1$, where
\beqano
\frac{1}{N_0}:= {
C_*\max
\left\{
\frac{\b_*}{c_0^2\varepsilon_0^2\d s_0}\sqrt{\frac{1}{\a_-}}\	 ,	\
\frac{\b_*}{c_0^2\varepsilon_0^{\frac{5}{2}}}{\frac{1}{\a_-}}
\right\}\frac{\a_+^{3/2}}{\a_-^{3/2}}
}
\eeqano
 with  $\beta_*$ as in~\equ{b*}
 and {$C_*$} is independent  of $s_0$.
Assuming also that \beqa{N0}N_0^{-1}<\frac{c_0^2\varepsilon_0^2\a_-^2}{2\a_+^2}\eeqa
we have, in particular, $N_0>2$.
We denote as

\beqa{stability}\widehat\HH_{*}=\hh(y_*)+g_{*}(y_*, x_*, \cG_*)+f_{*}(\cG_*, \g_*, y_*, x_*)\eeqa
the Hamiltonian obtained after the application of Theorem~\ref{NFL}
where, $g_*$, $f_*$ satisfy the  following bounds:
\beqano\|g_*{-\ovl f}\|\le 2 \D\, , \quad \|g_*\|\le 2^{-N}\D \eeqano
with $\ovl f(y_*, x_*, \cG_*)$ the $\g_*$--average of $f(y_*, x_*, \cG_*, \g_*)$ and 
\beqano
\D:=C_*\frac{m_0^2a\b_*}{c^2_0\varepsilon^2_0\a^2_-}
\eeqano
is an upper bound to $\|f\|$ above.
Let now $\G_*(t)=(\cG_*(t), \g_*(t), y_*(t), x_*(t)) $ be a solution of $\widehat\HH_*$ with initial datum $\G_*(0)$ $=$ $(\cG_*(0)$, $\g_*(0)$, $y_*(0)$, $x_*(0))$ $ \in{\mathbb D}$ and verifying
\beqa{initial data}&&|\cG_*(0)-1|\le \frac{\d}{2}\, , \quad 2\sqrt{m_0^3\a_-}\le|y_*(0)|{\le\frac{\sqrt{m_0^3\a_-}+\sqrt{m_0^3\a_+}}{2} }\nonumber\\
&& x_*(0)=\p\,.
\eeqa
We look for a time $T>0$ such that $\G_*(t)\in{\mathbb D}$
for all  $0\le t\le T$.  
Then $T$ can be taken to be
	\beqa{TOLD}T=\min\left\{\sqrt{{\a_-^3}}\, , \ {\frac{\sqrt{\a_-\varepsilon_0}}{\D}}\, , \ 2^{N_0}\frac{s_0\d}{\D}\right\}\eeqa
	as this choice easily allows to check
	\beqano
	&& |y_*(t)-y_*(0)|\le |y_*(0)|-\sqrt{\a_-}\,,\quad {|}\cG_*(t)-\cG_*(0)|\le \frac{2^{-(N+1)}\D t}{s_0}\le \frac{\d}{2}\nonumber\\
	&& |x_*(t)-x_*(0)|\le \p-\sqrt{\varepsilon_0}\eeqano
	for all $|t|\le T$. In addition, at the time $t=T$, one has
\beqano
|\g_*(T)-\g_*(0)|&\ge& c^\circ \min\left\{\b_*\sqrt{\frac{\a_-^3 }{\a_+^4}}\, , \ \frac{c^2_0\varepsilon^{5/2}_0\a^2_-}{\a_+^2}\sqrt{{\a_-}}\, , \ 
\frac{c^2_0\varepsilon^{2}_0\a^2_-}{\a_+^2}2^{N_0} s_0\d
\right\}\nonumber\\
&=:&\frac{3\p}{\eta}
\eeqano
with  $c^\circ$ independent of $\a_-$, $\a_+$, $\b$, $\ovl\b$, $\d$, $\varepsilon_0$ and $s_0$.
We then see that   $|\g_*(T)-\g_*(0)|$ is lower bounded by $3\p$ as soon as 
\beqa{eta}
\eta<1\,.
\eeqa 
The last step of the proof consists of proving that inequalities~\equ{ass1},~\equ{N0} and~\equ{eta} may be simultaneously satisfied. This is discussed in~\cite[Remark 5.1]{pinzari20a}. $\qquad \square$

\section{Chaos in a binary asteroid system}\label{Chaos in a binary asteroid system}
{This section describes the main steps of~\cite{diruzzaDP20} in constructing explicitly a topological horseshoe; henceforth providing evidences of the existence of symbolic dynamics.  
The construction essentially relies on the introduction of a two-dimensional Poincar\'e map from which invariants are computed. Following arguments presented in~\cite{ZGLICZYNSKI200432}, the introduction of ad hoc sets and the computation of their images provide the self--covering relationships needed to conclude.\\}

\noindent We fix $\b=\ovl\b=80$, which corresponds to take $\m=1$ and $\k\sim 40$; see~\equ{betagamma1}. We interprete the Hamiltonian $\widehat\HH_\JJ$
 with this choice of parameters as governing the (averaged out after many periods of the reference asteroid) motions  of a binary asteroid system interacting with a massive body, with the Jacobi reduction referred at one of the two twin asteroids. The parity triggered  by the equality $\b=\ovl\b$ reflects on the Taylor--Fourier coefficients of the expansion
 
 \begin{eqnarray}\label{hatHJ} 
\widehat\HH_\JJ (\RR, \GG, \rr, {\rm g})  
			 = \frac{\RR^2}{2}+ \frac{(\CC-\GG)^2}{2\rr^2}-\frac{1}{\rr}
+
\su\rr\sum_{\nu=1}^{\infty} q_{\nu}(\GG, {\rm g}) \left( \frac{\b}{\rr}\right)^{\nu}
\end{eqnarray}
accordingly to 
\begin{eqnarray}\label{qnu}
q_{\nu}(\GG, {\rm g}) =\left\{
\begin{aligned}
	&\sum_{p=0}^{\nu/2} \tilde q_{p}(\GG) \cos 2p\,{\rm g} \quad&{\rm if}\ \nu\ {\rm is\ even} \,   %\notag 
	\\
	 &\ 0 \quad &{\rm otherwise}\, . 
\end{aligned}
\right.	
\end{eqnarray}

\nl
In our numeric implementations, we  truncated the infinite sum in~\equ{hatHJ} up to a certain order  $\n_{\rm max}$,  chosen so that
the results did not vary increasing it again. Moreover, as  the coefficients $q_{\nu}(\GG,\cdot)$ in~\equ{qnu} are $\p$--periodic,  without loss of generality, we restricted ${\rm g}\in \torus/2\sim [0, \p)$. We describe the steps we followed in our numerical analysis,  recalling the reader to use Footnote~\ref{comparison} to relate the values in~\equ{starredvalues} and~\equ{q1q2} with the homonymous ones in~\cite{diruzzaDP20}.

\paragraph{Construction of a 2--dimensional Poincar\'e map}
The motions of~\equ{hatHJ} evolve on $3$--dimensional
manifolds $\cM_c$ labeled by the constant value $c$ of the energy. The structure of $\widehat\HH_\JJ$ allows to reduce the coordinate $\RR$ after fixing $c$ and hence to identify $\cM_c$ as the 3--dimensional space of triples $\{(\rr, \GG, {\rm g})\}$. The dimension can be further reduced to $2$ considering a plane $\Sigma$ through a given
$P_*=(\rr_*, \GG_*, {\rm g}_*)$ and perpendicular to the velocity vector $V_*=(v^*_\rr, v^*_\GG, v^*_{\rm g})$ of the the orbit through $P_*$. This leads us to construct 
a Poincar\'e map, which we define as follows. We start by defining two operators $l$ and $\pi$ consisting in ``lifting'' the initial two--dimensional seed $z=(\GG,{\rm g})$ to the four--dimensional space $(\RR, \GG, \rr, {\rm g})$ and ``projecting'' it back to plan after the action of the flow--map $\Phi^{t}_{\widehat\HH_\JJ}$ during the first return time $\tau$. The lift operator reconstructs the four--dimensional state vector from a seed on $D \times \mathbb{T}/2$, where the domain $D$ of the variable $\GG$ is a compact subset of the form $[-1,1]$. For a suitable $(\mathcal{A},A) \subset \R^{2} \times \R^{2} $, its definition reads 
\begin{eqnarray}
l: \quad  D \times \mathbb{T}/2 \supset \mathcal{A}   & \rightarrow & D  \times \mathbb{T}/2 \times A \notag \\
   	      z	     & \mapsto &  \tilde{z}=l(z) \notag,
\end{eqnarray}
where  $\tilde{z}=(\GG,{\rm g},\RR,\rr)$ satisfies the two following conditions:
\begin{enumerate}
	\item \textit{Planarity condition.} The triplet $(\rr, \GG,{\rm g})$ belongs to the plane $\Sigma$, \ie $\rr$ solves the algebraic condition 
	$ v_{\rr}^{*} (\rr-\rr_{*})+v_{\GG}^{*} (\GG-\GG_{*})+ v_{\rm g}^{*} (\rm g-\rm g_{*}) = 0$.  
	\item \textit{Energetic condition.} The component $\RR$ solves the energetic condition $\widehat\HH_\JJ(\RR, \GG, \rr, {\rm g})=c$.  
\end{enumerate} 
The projector $\pi$ is the projection onto the first two components of the vector,
\begin{eqnarray*}
\pi: \quad  D \times \mathbb{T}/2 \times A &\rightarrow & D \times \mathbb{T}/2  \notag \\
   	        \tilde{z}=(z_{1},z_{2},z_{3},z_{4})	     &\mapsto  & \pi(z)=(z_{1},z_{2})\, .
\end{eqnarray*}
The Poincar\'e mapping is  therefore defined and constructed as 
\begin{eqnarray}
P: \quad D \times \mathbb{T}/2 &\rightarrow & D \times \mathbb{T}/2 \notag \\
   	      z	     &\mapsto &  z'=P(z)= 
	      \big(\pi \circ \Phi^{\tau(z)}_{\widehat\HH_\JJ} \circ l \big)
	      (z)\, .\notag
\end{eqnarray}

\nl
 The mapping is nothing else than a ``snapshots'' of the whole flow at specific return time $\tau$. It should be noted that the successive (first) return time is in general function of the current seed (initial condition or current state), \ie $\tau=\tau(\tilde{z})$, formally defined (if it exists) as 
\begin{eqnarray*}
	\tau(z) = \inf \Big\{ t \in \mathbb{R}_{+}, 
	\big(\rr(t), \GG(t), {\rm g}(t)\big)	\in \Sigma   \Big\},
\end{eqnarray*}
where $\big(\rr(t), \GG(t), {\rm g}(t)\big)$ is obtained through the flow. 

\begin{figure}[t!]
	\centering
	\includegraphics[width=0.95\linewidth]{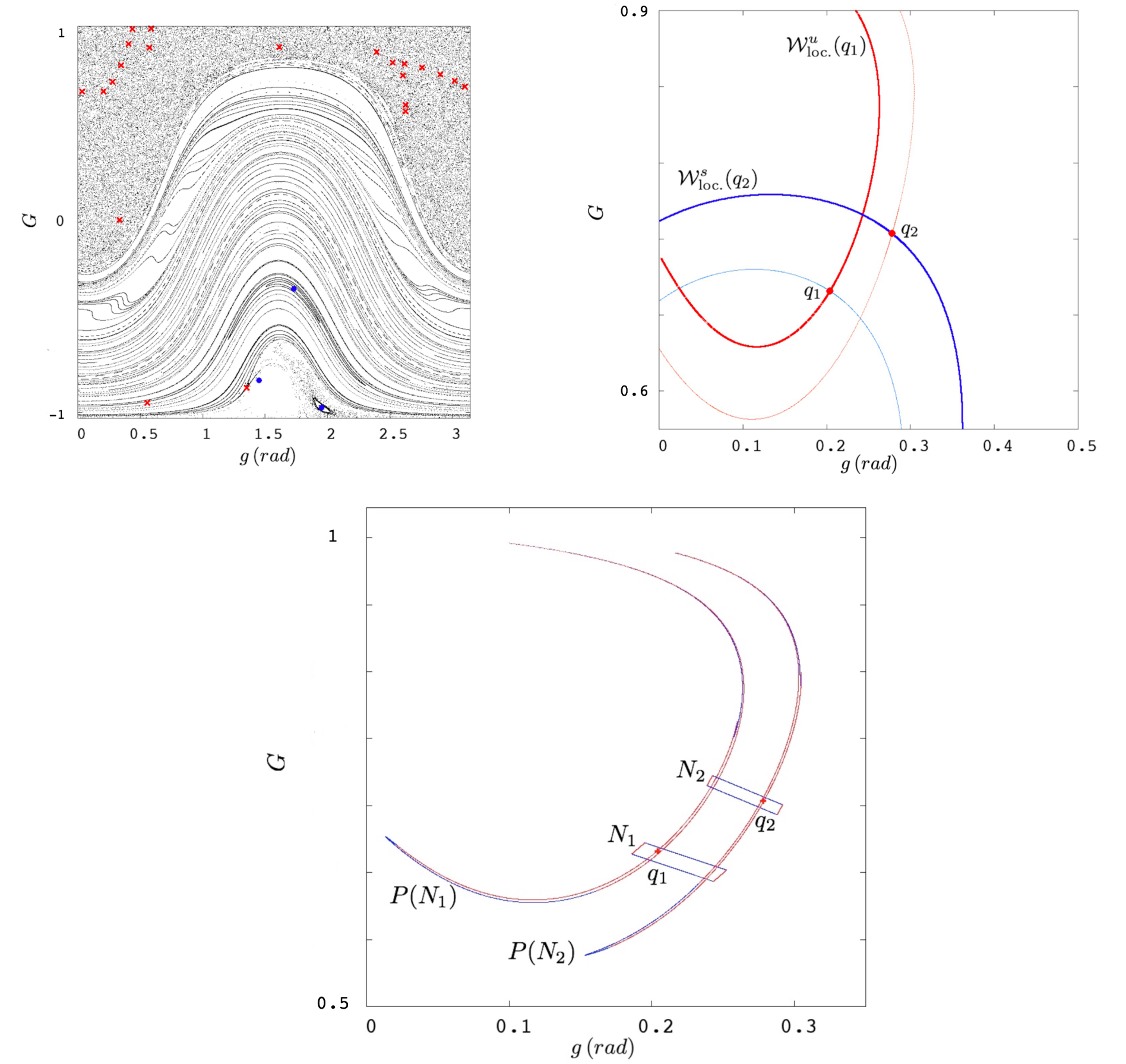}
	\caption{{Composite panels illustrating our main steps in constructing the topological horseshoe. (Top left) The continuous flow is reduced to a 2--dimensional mapping by introducing a suitable Poincar\'e map $P$. The phase space contains both elliptic (blue) and hyperbolic (red) fixed-points. (Top right) Finite pieces of the stable and unstable manifolds might be constructed from the knowledge of the eigensystem derived from the linearisation $DP$. (Bottom) Carefully chosen sets and their images under $P$ provide the covering relations and imply existence of symbolic dynamics. See text for details.}} 
	\label{fig:fig3}
\end{figure}

	\nl
	With $\CC=24.394$, we fixed the initial values
\begin{eqnarray}\label{starredvalues}
 \RR_{*} = %\cancel{-0.0039}
 {-0.0060}\,,\quad \GG_{*} = %\cancel{-2.4915}
 {-0.804}\,, \quad  \rr_{*}  = %\cancel{3132.069}
 {652.256}\,,\quad {\rm g}_{*}  =  1.4524\ \rm rad
\end{eqnarray}
and we obtained the results plotted in {the first panel of} Figure~\ref{fig:fig3}. We invite the reader to compare this figure with the unperturbed phase portrait of Figure~\ref{fig:fig2}.
In particular, due to the non--integrability of the problem,  chaotic zones appear, mostly distributed   {for positive values of $G$}.
This chaos was the object of our next investigations, as discussed below.

\paragraph{Hyperbolic fixed points and heteroclinic intersections}
Equilibrium points of the mapping $P$ (\ie periodic orbits of the Hamiltonian system~\equ{hatHJ}) have been found using a Newton algorithm with initial guesses distributed on a resolved grid of initial conditions in $D \times \mathbb{T}/2$.  We found more than $20$ fixed points $x_{*}$.  The eigensystems associated to the fixed points have been computed to determine the local stability properties. 
 The result of the analysis is displayed on Figure~\ref{fig:fig3} along with the following convention: hyperbolic fixed points appear as red crosses, elliptical points are marked with blue circles. \\
  The local stable manifold associated to an hyperbolic point $x_{\star}$,
\begin{eqnarray*}
\mathcal{W}^{s}_{\textrm{loc.}}(x_{\star}) = 
	\Big\{ 
	x \, \vert \,
	\|{P^{n}(x)-x_{\star}}\| \to 0, \, n \in \mathbb{N}_{+}, n \to \infty
	\Big\},
\end{eqnarray*}
can be grown by computing the images of a fundamental domain $I \subset E_{s}(x_{\star})$,  $E_{s}(x_{\star})$ being the stable eigenspace associated to the saddle point $x_{\star}$. 
 The local unstable manifolds $\mathcal{W}^{u}_{\textrm{loc.}}(x_{\star})$ were similarly computed, but changing the sign of the time integration. See Figure~\ref{fig:fig3}.

\paragraph{Covering relations}

Let us introduce some notations. Let $N$ be a compact set contained in 
$\R^2$ and $u(N)=s(N)=1$ being, respectively, the {\it exit} and {\it entry dimension} (two real numbers such that their sum is equal to the dimension of the space containing $N$); let $c_{N} : \R^2 \rightarrow \R^2$ be an homeomorphism such that $c_N(N) = [-1,1]^2$; let $N_c=[-1,1]^2$, $N_c^-=\{-1,1\} \times [-1,1]$, $N_c^+= [-1,1]\times \{-1,1\}$; then, the two set $N^-= c_N^{-1}(N_c^-)$ and $N^+= c_N^{-1}(N_c^+)$ are, respectively, the {\it exit set} and the {\it entry set}. In the case of dimension $2$, they are topologically a sum of two disjoint intervals.  The quadruple $(N,u(N),s(N),c_N)$ is called a {\it h--set} and $N$ is called {\it support} of the $h$--set. Finally, let $S(N)_c^l =(-\infty, -1) \times \R$, \,  $S(N)_c^r =(1,\infty) \times \R$, and $S(N)^l = c_N^{-1}(S(N)_c^l) , \, S(N)^r = c_N^{-1}(S(N)_c^r) $ be, respectively, the left and the right side  of $N$.   The general definition of covering relation can be found in~\cite{aGi19}. Here we provide a simplified notion, suited to the case that $N$ is two--dimensional, based on~\cite[Theorem 16]{ZGLICZYNSKI200432}.
\begin{definition}\rm
\label{def_covering}
	Let $f : \R^2 \rightarrow \R^2$ be a continuous map and $N$ and $M$ the supports of two $h$--sets. We say that $M$ $f$--covers $N$ and we denote it by $M \stackrel{f} \Longrightarrow N $ if:
	\begin{itemize}
	\item[(1)] $\exists\, q_0\in [-1, 1]$ such that $f(c_N([-1, 1]\times \{q_0\}))\subset {\rm int}( S(N)^l \bigcup N \bigcup S(N)^r)$,
	\item[(2)] $f(M^-) \bigcap N  = \emptyset$,
	\item[(3)] $f(M) \bigcap N^+  = \emptyset$.
\end{itemize}	
Conditions (2) and (3) are called, respectively, {\it exit} and {\it entry condition}.
\end{definition}\rm
\noindent  The notions of covering (including self--covering) relations are useful in defining \textit{topological horseshoe}~\cite{aGi19},~\cite{ZGLICZYNSKI200432}.\\

\begin{definition}\rm
	Let $N_{1}$ and $N_{2}$ be the supports of two disjoint $h$--sets in $\R^2$. A continuous map $f : \R^2 \rightarrow \R^2$ is said to be a {\it topological horseshoe} for $N_1$ and $N_2$ if 
	\begin{eqnarray*}
		N_1 \stackrel{f} \Longrightarrow N_1 \, , \quad N_1 \stackrel{f} \Longrightarrow N_2 \, , \quad
		N_2 \stackrel{f} \Longrightarrow N_1 \, , \quad N_2 \stackrel{f} \Longrightarrow N_2 \, .
	\end{eqnarray*}
\end{definition}
\noindent Topological horseshoes are associated to symbolic dynamics as discussed in 
in~\cite[Theorem 2]{aGi19} and in~\cite[Theorem 18]{ZGLICZYNSKI200432}.

\paragraph{The topological horseshoe}

 Based on the couple of hyperbolic fixed points 
\begin{eqnarray}\label{q1q2}
\left\{
\begin{aligned}
&q_{1}=({\rm g}_{1},\GG_{1}) = (0.203945459,%\cancel{2.06302430} 
{0.665706}), \\ 
&q_{2}=({\rm g}_{2},\GG_{2}) = (0.278077917,%\cancel{2.21418596} 
{0.714484}),
\end{aligned}
\right.
\end{eqnarray}
we define  two sets $N_{1}, N_{2} \subset \R^2$ which are supports of two $h$--sets as follows:
\begin{eqnarray}
\left\{
\begin{aligned}
&N_{1} = q_1 + A_1 v_1^s +B_1 v_1^u, \\ \notag 
&N_{2} = q_2 + A_2 v_2^s +B_2 v_2^u,
\end{aligned}
\right.
\end{eqnarray}
where 
 $v_{1}^s$, $v_1^u$, $v_2^s$, $v_2^u$ are the stable and the unstable eigenvectors related to $q_{1}$, $q_{2}$, respectively, and the $A_1$, $A_2$,  $B_1$ and $B_2$ are numbers suitably chosen in a grid of values. Then the following covering relations are numerically detected
\begin{eqnarray*}
	N_1 \stackrel{P} \Longrightarrow N_1 \, , \quad N_1 \stackrel{P}\Longrightarrow N_2 \, , \quad
	N_2 \stackrel{P} \Longrightarrow N_1 \, , \quad N_2 \stackrel{P} \Longrightarrow N_2  \, ,
\end{eqnarray*}
proving the numerical evidence  of a topological horseshoe, \ie, existence of symbolic dynamics for $P$. The obtained horseshoe associated to $q_{1}$ and $q_{2}$ with the aforementioned parameters is illustrated in Figure\,\ref{fig:fig3}.

\appendix
\section{Outline of the proof of Theorem~\ref{NFL}}\label{Outline of the proof}

In this section we {provide some technical details of} the proof of Theorem~\ref{NFL}. For the full proof we refer to~\cite{pinzari20a}.

\nl
The proof is by recursion.
We assume that, at a certain step, we have a system of the form
\beqa{step i}\HH({\mathbf I}, \bm\f, {\mathbf J}({\mathbf p}, {\mathbf q}), y)=\hh({\mathbf I}, {\mathbf J}({\mathbf p}, {\mathbf q}), y)+g({\mathbf I}, {\mathbf J}({\mathbf p}, {\mathbf q}), y, x)+f({\mathbf I}, \bm\f, {\mathbf J}({\mathbf p}, {\mathbf q}), y, x)\eeqa
where  $f\in {\cal O}_{\r, s, \d, r, \xi}$, $g\in {\cal N}_{\r, s, \d, r, \xi}$. At the first step, just take $g\equiv 0$.

\nl
After splitting $f$ on its Taylor--Fourier basis 
\beqano f=\sum_{k,h,j} f_{khj}({\mathbf I},y, x)e^{\ii k \cdot\bm{\f}}{\mathbf p}^h  {\mathbf q}^j\, . \eeqano
one looks for a time--1 map 
\beqano\Phi=e^{\cL_\phi}=\sum_{k=0}^{\infty}\frac{\cL_\phi^k}{k!}\qquad \cL_\phi(f):=\big\{\phi,\ f\big\}\eeqano generated by a small Hamiltonian 
$\phi$
which will be taken  in the  class ${\cal Z}_{\r, s, \d, r, \xi}$ in~\equ{zero average}.
Here,
\beqano \big\{\phi,\ f\big\}&:=&\sum_{i=1}^n(\partial_{{\mathbf I}_i}\phi\partial_{\bm\f_i}f-\partial_{{\mathbf I}_i}f\partial_{\bm\f_i}\phi)+\sum_{i=1}^m(\partial_{{\mathbf p}_i}\phi\partial_{\mathbf q_i}f-\partial_{{\mathbf p}_i}f\partial_{\mathbf q_i}\phi)\nonumber\\
&&+(\partial_{y}\phi\partial_{x}f-\partial_{y}f\partial_{x}\phi)\eeqano
denotes the Poisson parentheses of $\phi$ and $f$.
One lets
\beqa{exp}\phi=\sum_{(k,h,j):\atop{(k,h-j)\ne (0,0)}} \phi_{khj}({\mathbf I},y, x)e^{\ii k\cdot\bm{\f}}{\mathbf p}^h {\mathbf q}^j\, . \eeqa

\nl
The operation
\beqano \phi\to \{\phi,\hh\}\eeqano
acts diagonally  on the monomials  in the expansion~\equ{exp}, carrying
\beqa{diagonal}\phi_{khj}\to -\big(\o_y\partial_x \phi_{khj}+\l_{khj} \phi_{khj}\big)\, , \quad {\rm with}\quad \l_{khj}:=(h-j)\cdot\o_{\mathbf J}+\ii k\cdot\o_{\mathbf I}\, . \eeqa
Therefore, one defines
\beqano\{\phi,\hh\}=:-D_\omega\phi\, . \eeqano
The formal application of $\Phi=e^{\cL_\phi}$ yields:
\beqa{f1}
e^{\cL_\phi} \HH&=&e^{\cL_\phi} (\hh+g+f)=\hh+g-D_\omega \phi+f+\Phi_2(\hh)+\Phi_1(g)+\Phi_1(f)
\eeqa
where the $\Phi_h:=\Phi_h:=\sum_{j\ge h}\frac{{\cal L}_\phi^j}{j!}$'s are the tails of $e^{\cL_\phi}$.

\nl
Next, one requires that   the residual term $-D_\omega \phi+f$ lies in the class ${\cal N}_{\r, s, \d, r, \xi}$ in~\equ{phi independent}%. This
%amounts  to solve 
%the ``homological'' equation
\beqa{homological equation}{(-D_\omega \phi+f
)}\in {\cal N}_{\r, s, \d, r, \xi}\eeqa
for $\phi$.

\nl
Since we have chosen $\phi\in{\cal Z}_{\r, s, \d, r, \xi}$, by~\equ{diagonal}, we have that also $D_\omega \phi\in{\cal Z}_{\r, s, \d, r, \xi}$. So, Equation~\equ{homological equation} becomes
\beqano-D_\omega \phi+\widetilde f=0\, . \eeqano
In terms of  the Taylor--Fourier modes, the equation becomes
\beqa{homol eq} \o_y\partial_x \phi_{khj}+\l_{khj} \phi_{khj}=f_{khj}\qquad \forall\ (k,h,j):\ (k,h-j)\ne (0,0)\, . \eeqa

\nl In the standard situation, one typically proceeds to solve such equation via Fourier series:
\beqano f_{khj}({\mathbf I},y, x)=\sum_{\ell}f_{khj\ell}({\mathbf I}, y)e^{\ii \ell x}\, , \qquad \phi_{khj}({\mathbf I},y, x)=\sum_{\ell}\phi_{khj\ell}({\mathbf I}, y)e^{\ii \ell x}\eeqano
so as to find
$\dst\phi_{khj\ell}=\frac{f_{khj\ell}}{\m_{khj\ell}}$
with the usual  denominators $\m_{khj\ell}:=\l_{khj}+\ii\ell \omega_y$ which one requires not to vanish %(equivalently, $\frac{\l_{khj}}{\o_y}\notin \ii\integer$)
via, e.g., a ``diophantine inequality'' to be held for all $(k,h,j,\ell)$ with $(k,h-j)\ne (0,0)$. In this standard case, there is not much freedom in the choice of $\phi$. In fact, such solution is determined up to solutions  of the homogenous equation
\beqa{homogeneous}D_\omega\phi_0=0\eeqa
which, in view of the Diophantine condition, has the only trivial solution $\phi_0\equiv0$. {\it The situation is different if $f$ is not periodic in $x$, or $\phi$ is not needed so}. In such a case, it is possible to find a solution of~\equ{homol eq}, corresponding to a non--trivial solution of~\equ{homogeneous},  where small divisors do not appear. This is
\beqa{solution}\phi_{khj}({\mathbf I}, y,  x)=\arr{\dst
\frac{1}{\o_y}\int_0^xf_{khj}({\mathbf I}, y, \t)e^{\frac{\l_{khj}}{\o_y}(\t-x)}d\t\quad {\rm if}\ \ (k,h-j)\ne (0,0)\\\\
\dst0\qquad\qquad\qquad\qquad\qquad\qquad\qquad\ {\rm otherwise .}
}
\eeqa
Multiplying by $e^{ik\f}$ and summing over $k$, $h$ and $j$, we obtain
\beqano\phi({\mathbf I}, {\bm \f}, p, q, y, x)=\frac{1}{\o_y}\int_0^x \widetilde f\left({\mathbf I}, {\bm \f}+\frac{\o_{\mathbf I}}{\o_y}(\t-x), p e^{\frac{\o_{\mathbf J}}{\o_y}(\t-x)}, q e^{-\frac{\o_{\mathbf J}}{\o_y}(\t-x)}, y, \t\right)d\t\, . \eeqano
In~\cite{pinzari20a} it is proved  that, under the assumptions~\equ{normal form assumptions}, this 
function can be used
to obtain a convergent time--one map and that the construction can be iterated so as to provide the proof of Theorem~\ref{NFL}. The construction of the iterations and the proof of its convergence is obtained adapting the techniques of
~\cite{poschel93} to the present case.
%{Manca il comando {\tt addcontentsline\{toc\}\{section\}\{References\}}}

\subsection*{Acknowledgments} The authors acknowledge the  European Research Council (Grant 677793 Stable and Chaotic Motions in the Planetary Problem) for supporting them during the completion of the results described in the paper; warmly thank the organisers of  I--CELMECH Training School that held in Milan in winter 2020 for their interest and especially U.~Locatelli for a highlighting discussion. 
\addcontentsline{toc}{section}{References}
	 \bibliographystyle{plain}
%\bibliography{REFERENCES.bib}\def\cprime{$'$} \def\cprime{$'$}

\def\cprime{$'$} \def\cprime{$'$}

\end{document}